\documentclass[10pt,aps,pre,noshowpacs,twocolumn,superscriptaddress,nobibnotes,nofootinbib,floatfix,a4paper]{revtex4}
\usepackage{amssymb}
\usepackage{graphicx}
\usepackage{amsmath,amsfonts}
\usepackage[colorlinks=true,allcolors=blue]{hyperref}
\usepackage[utf8]{inputenc}
\usepackage{xcolor}  \usepackage[normalem]{ulem}
\usepackage[T1]{fontenc} 
\usepackage{multirow}
\usepackage{booktabs}
\usepackage{threeparttable}
\usepackage{tabularx}
\usepackage{subcaption}

\usepackage{pdfpages}

\usepackage{appendix}

\usepackage{mathtools}
\usepackage[normalem]{ulem}
\usepackage{url}
\usepackage[colorlinks=true,allcolors=blue]{hyperref}
\usepackage{float}
\usepackage{wasysym}
\usepackage{bbm}

\usepackage{booktabs}   
\usepackage{array}      

\definecolor{rj}{HTML}{FF9966}

\usepackage{placeins}

\begin{document}

\title{Motif-based filtrations for persistent homology: \\
A framework for graph isomorphism and property prediction}

\author{Meritxell Vila-Mi\~nana}
\affiliation{Center for Complex Networks and Systems Research, Luddy School of Informatics, Computing, and Engineering, Indiana University, Bloomington, IN, USA}
\author{Robert Jankowski}
\affiliation{Faculty of Electrical Engineering, Mathematics and Computer Science, Delft University of Technology, 2628 CD, Delft, Netherlands}
\affiliation{Departament de F\'isica de la Mat\`eria Condensada, Universitat de Barcelona, Mart\'i i Franqu\`es 1, E-08028 Barcelona, Spain}
\affiliation{Institute of Complex Systems (UBICS), Universitat de Barcelona, Barcelona, Spain}
\author{Aina Ferr\`a Marc\'us}
\affiliation{Departament de Mat\`ematiques i Inform\`atica, Universitat de Barcelona,Gran Via de les Corts Catalanes 585, E-08007 Barcelona, Spain}
\author{Rub\'en Ballester}
\affiliation{Departament de Mat\`ematiques i Inform\`atica, Universitat de Barcelona,Gran Via de les Corts Catalanes 585, E-08007 Barcelona, Spain}
\author{M. {\'A}ngeles Serrano}
\affiliation{Departament de F\'isica de la Mat\`eria Condensada, Universitat de Barcelona, Mart\'i i Franqu\`es 1, E-08028 Barcelona, Spain}
\affiliation{Institute of Complex Systems (UBICS), Universitat de Barcelona, Barcelona, Spain}
\affiliation{ICREA, Passeig Llu\'is Companys 23, E-08010 Barcelona, Spain}
\author{Carles Casacuberta}
\affiliation{Departament de Mat\`ematiques i Inform\`atica, Universitat de Barcelona,Gran Via de les Corts Catalanes 585, E-08007 Barcelona, Spain}


\begin{abstract}
Determining whether two graphs are isomorphic is a fundamental problem with practical applications in areas such as molecular chemistry or social network analysis, yet it remains a challenging task, with exact solutions often being computationally expensive. We address this task using persistent homology built on motif-based filtrations of graphs, a method from topological data analysis that summarizes the shape of data by tracking the persistence of structural features along filtrations. Specifically, we use edge-weighting schemes based on the densities of triangles, chordless squares, and chordless pentagons, which have been shown to be effective for detecting network dimensionality. Our cycle-density filtrations distinguish non-isomorphic graphs perfectly or nearly perfectly across four demanding graph families, many of which exhibit symmetries. We outperform curvature-based, degree-based, and Vietoris--Rips filtrations, and match or exceed the accuracy of egonet-distance methods while incurring a lower computational cost. The expressive power of our filtrations goes beyond isomorphism testing: because they capture rich structural information from graphs, they consistently achieve top performance on property prediction tasks using real-world data, and exhibit high sensitivity to edge rewiring and removal. Together, these findings establish cycle-density filtrations as an effective and computationally tractable framework for graph comparison and characterization, bridging topological data analysis and network science.
\end{abstract}

\maketitle

\let\oldaddcontentsline\addcontentsline
\renewcommand{\addcontentsline}[3]{}

\section{Introduction}
\label{sec:introduction}

Distinguishing whether two graphs are isomorphic is a fundamental problem with applications in different scientific domains. For instance, the verification of molecular structures and chemical compound descriptions emphasizes the need for efficient and expressive algorithms to test graph isomorphism \cite{akutsu2013comparison,merkys2023graph}. 
In biomedical research, graph isomorphism has been applied to analyze protein structures, and more broadly in social and biological networks for pattern recognition \cite{sangkaran2019survey}.

Recently, topology-based methods have emerged as a promising direction for graph comparison and isomorphism testing. In particular, persistent homology, a central tool in topological data analysis (TDA), provides descriptors capable of capturing structural properties of graphs. The effectiveness of persistent homology depends crucially on the choice of filtration, which determines how a graph is transformed into a sequence of nested graphs or higher-dimensional simplicial complexes that expose multiscale topological features of the underlying data. Despite growing interest in TDA, the impact of filtration design on the performance of several tasks, such as distinguishing non-isomorphic graphs, graph learning, or robustness to small perturbations, remains insufficiently understood~\cite{Ballester2024}. 

In this work, we use persistent homology computed from filtrations based on the densities of chordless cycles to evaluate graph isomorphism and predict graph properties. 
Edge weightings based on chordless-cycle densities (triangles, squares, and pentagons) were introduced in~\cite{almagro2022detecting} for latent dimensionality detection in complex networks, and were integrated with TDA techniques in \cite{ferra2025} to further strengthen dimensionality detection capabilities. 
Densities of chordless cycles are related to Ollivier--Ricci curvature \cite{ollivier2007ricci}, which is a discrete analog of Ricci curvature in Riemannian geometry. Positive curvature in a graph region is associated with a high density of triangles, while negative curvature is linked to increased densities of chordless longer cycles~\cite{lin2011ricci,Jost-Liu_2014}. 
Curvature defines local geometric invariants that are preserved under graph isometries, and has been shown to be effective for identifying community boundaries, i.e., edges separating densely connected regions, as well as for assessing the robustness and vulnerability of networks~\cite{Sreejith_etal_2016,Sandhu_etal_2015,Jost-Liu_2014,Fesser_2024}.

Our results demonstrate that persistence descriptors computed from chordless-cycle densities are highly accurate in distinguishing non-isomorphic graphs, particularly in challenging datasets. In doing so, we extend the framework introduced in~\cite{Ballester2024} to study the expressivity of persistent homology for graph learning by comparing several filtrations inspired by graph geometry and curvature. 
In addition to the filtrations considered in \cite{Ballester2024} (a degree-based filtration, 
filtrations based on the Ollivier--Ricci and augmented Forman--Ricci curvatures, and a Vietoris--Rips filtration defined by the shortest-path distance between nodes) we evaluate edge-based filtrations derived from density of triangles, density of chordless squares, density of chordless pentagons, and sum of densities of triangles, squares and pentagons \cite{almagro2022detecting}, as well as classical network measures such as the Randi\'{c} connectivity index \cite{randic1975characterization}, weights used to compute the harmonic index \cite{fajtlowicz1987conjectures} and the repulsion-attraction rule \cite{muscoloni2017machine}, betweenness centrality \cite{muscoloni2017machine}, clustering coefficient \cite{watts1998collective}, egonet persistence \cite{piccardi2023metrics}, 
and a graphlet-based filtration \cite{sarajlic2016graphlet}. 

Beyond the graph isomorphism problem, we evaluate the performance of filtrations on real-world datasets using property prediction tasks, where the goal is to predict structural characteristics of graphs directly from their persistent homology features. In particular, we consider global and local network descriptors, including average clustering coefficient, average shortest-path length, degree statistics, and centrality measures. This setting allows us to assess how well persistent homology captures meaningful structural information beyond its ability to distinguish non-isomorphic graphs. We also analyze the sensitivity of filtrations to small perturbations of graphs, such as rewiring and edge removal. 

Across benchmark datasets, our motif-based filtrations consistently outperform filtrations based on degree, curvature, and centrality measures in distinguishing non-isomorphic graphs, including strongly regular graphs where standard methods fail. Similarly, our filtrations achieve the highest performance in property prediction and sensitivity analyses, using both synthetic and real data.

\subsection{Related Work}
Graph isomorphism is closely related to the broader problem of network comparison, where the goal is quantifying similarity between graphs. In \cite{tantardini2019comparing} and \cite{hartle2020network}, several network comparison measures were studied. While checking whether two finite graphs are isomorphic is a complex task, they addressed the question on whether two real-world networks are similar.

Comparing networks requires an inexact graph matching framework, namely the definition of a real-valued distance that converges to zero as the networks become increasingly close to being isomorphic.
A wide variety of network comparison methods and algorithms were proposed in \cite{hartle2020network}, which presented a review of several approaches to network comparison. Among these, the best state-of-the-art techniques involve graphlet-based methods. Graphlets are small connected induced subgraphs, typically restricted to at most five nodes to limit computational complexity. Advanced graphlet-based methods account for automorphism orbits, thereby distinguishing the structural roles of nodes within each graphlet. Empirical evaluations in \cite{tantardini2019comparing} show that, for undirected networks, the Graphlet Correlation Distance GCD-11 \cite{yaverouglu2014revealing} is the best-performing distance in discriminating between different network topologies, outperforming all other methods, while for directed networks, DGCD-129 \cite{sarajlic2016graphlet} led to an optimal performance.

An alternative line of work focuses on egonet-based representations, which summarize local network structure around each node. In \cite{piccardi2023metrics}, Piccardi proposed an alignment-free network comparison method based on distributions of egonet features, including normalized degree, clustering coefficient, and egonet persistence. He suggested taking the statistics of these three indicators in different combinations to define a distribution
function in 1, 2, or 3 dimensions, as a synthesis of network properties. Similarly to graphlet-based techniques, the local graph structure around each node was summarized by a vector of features, which was used as a global descriptor of the network. 
Yet, the features used in \cite{piccardi2023metrics} are simpler and enable rapid computation. The results of their experiments show that ego-distances perform comparably to the best graphlet-based distance (GCD-11) with similar computational requirements.

A central tool for graph isomorphism detection is the Weisfeiler--Leman algorithm \cite{weisfeiler1968reduction}, 
whose use in a machine-learning setting was discussed in~\cite{morris2019weisfeiler,molina2025topological}.
From a topological perspective, Ballester and Rieck~\cite{Ballester2024} provided a theoretical study and empirical analysis of the expressivity of persistent homology for graph learning tasks. They evaluated the expressivity of five filtrations by letting them distinguish non-isomorphic graphs. In addition, they evaluated such filtrations for graph property prediction (e.g., diameter and girth) and graph classification tasks, demonstrating the utility of topological features on benchmark datasets.

Other studies have combined classical network measures with persistent homology to obtain a scalable graph descriptor and used it to compare similarities between graphs. For example, in \cite{hajij2020fast}, the PageRank algorithm was used along with persistent homology for graph similarity analysis, thus showing the usefulness of network measures for graph comparison. Similarly, a topological framework based on scale-space models was used in~\cite{molina2025topological} to analyze hypergraphs. Their approach uses hypergraph representations and $s^2$\nobreakdash-models, which involve a dynamic sequence of chain complexes and transition maps. 

A~triangle-aware filtration for persistent homology was used in \cite{calissano2025topological} to preserve triangular cycles in spatial graphs during topological coarsening. While conceptually related to our approach, their method relies on a modified clique filtration that delays the appearance of $2$\nobreakdash-simplices to repair the fact that triangles become irrelevant in clique filtrations. 
Persistent homology has also been applied to property prediction, where representations based on persistence diagrams are used as features for supervised learning~\cite{rieck2018neural, carriere2020perslay, tola2024toper}, highlighting the potential of topological descriptors for graph learning tasks when applied with suitable filtrations. 

\section{Methods}

\subsection{Persistent homology}
Persistent homology is a tool from algebraic topology that captures how structural features of a graph ---or, more generally, of a simplicial complex--- appear and disappear as a filtration parameter evolves.
In this work, persistent homology is used to summarize and compare information induced by different graph filtrations. 
A graph filtration is a sequence of nested subgraphs exhausting the whole graph, usually sublevel sets of a real-valued function on nodes and edges.
Persistent homology tracks the birth and death of topological features along a filtration, such as connected components and edge cycles. 
Cycles disappear (as homology generators) when higher simplices from the clique complex fill them, while connected components arise and merge along the filtration. 
This information is summarized in a persistence diagram, which consists of points $(b,d)$ representing the filtration values at which a feature appears (birth) and disappears (death). Features that persist over a wide range of scales are considered structurally significant, while short-lived features are often interpreted as noise.

A graph is denoted as $G=(N,E)$, where $N$ is its set of nodes and $E$ is its set of edges.
The graph isomorphism problem aims to determine whether, given two graphs $G=(N,E)$ and $G'=(N',E')$, there is a bijective correspondence $f\colon N\to N'$ such that $\{f(u),f(v)\}\in E'$ if and only if $\{u,v\}\in E$.
In our context, filtrations are equivariant in the sense of \cite{Ballester2024}, that is, filtering functions are compatible with graph isomorphisms and hence yield identical persistence diagrams on isomorphic graphs.

\begin{figure}[htbp]
    \centering
\includegraphics[width=0.88\columnwidth]{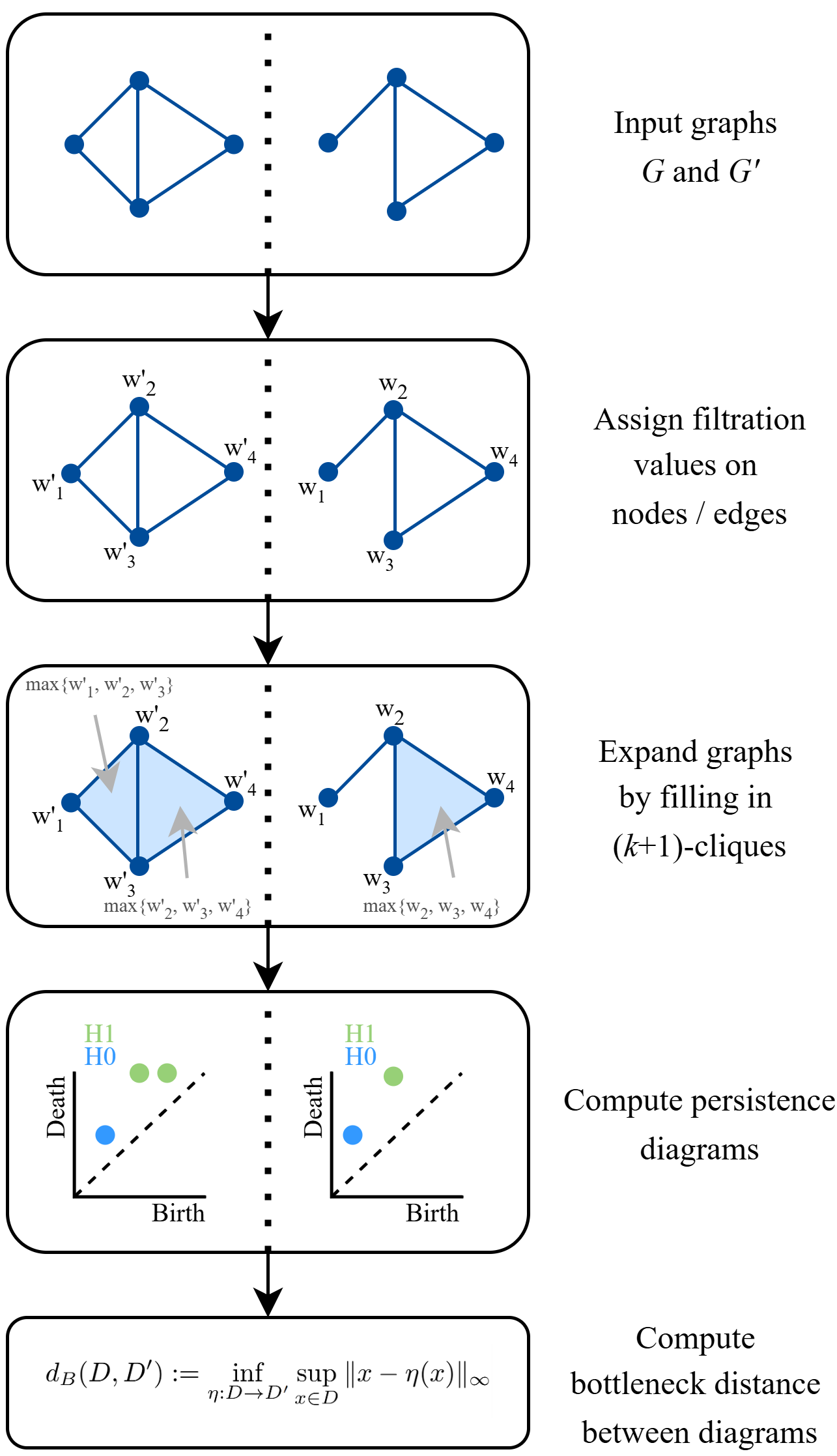}
\caption{Pipeline for the graph isomorphism problem.}
\label{fig:pipeline}
\end{figure}

Our approach to the graph isomorphism problem is depicted in Fig.\;\ref{fig:pipeline} and begins with taking two input graphs $G$ and $G'$.
The pipeline begins by assigning weights to the nodes or edges of $G$ and~$G'$ that serve as filtration values. If an edge weighting is given, then nodes are weighted with the absolute minimum weight of all the edges. If a node weighting is given instead, then each edge is weighted with the maximum value of its endpoints. Moreover, the graph is expanded into a simplicial complex by filling in its cliques (i.e., complete subgraphs), where each simplex is weighted with the maximum value of its faces.
Next, we compute a persistence diagram for the resulting simplicial complex with respect to the weight filtration, by considering the sequence of sublevel simplicial complexes.
Finally, we compute the bottleneck distance $d_B(D,D')$ between the corresponding persistence diagrams $D$ and~$D'$, and conclude that $G$ and $G'$ are isomorphic if $d_B(D,D') < 10^{-8}$.  
The full setup is described in more detail in the Appendix (Section\!~\ref{subsec:experimental_setup}).

One of the crucial parts of the algorithm is the weighting scheme of the graph, since different filtrations yield distinct
persistence diagrams. Thus, this article focuses on comparing various filtrations and analyzing which one is most effective 
for detecting non-isomorphic graphs, predicting graph properties, and assessing sensitivity to small changes. Table~\ref{tab:filtration_legend} provides a legend of the filtrations used in this paper, whose definitions are detailed in Section\!~\ref{subsec:def_filtrations} in the Appendix. These filtrations comprise several categories: we study local node-level scalar measures (nD, nC, nE), curvature-based measures (eO, eF), motif-based filtrations (eT, eS, eP, e$\Sigma$), distance-based measures (mV), centrality-based measures (eB), degree-based connectivity indices (eR, eH, eA), and graphlet-based filtrations (nG).

\begin{table}[htbp]
\renewcommand{\arraystretch}{1.3}
\caption{List of filtrations used in this work.}
\label{tab:filtration_legend}
\centering
\begin{tabular}{ll}
\hline
\multicolumn{2}{c}{\textbf{Edge-based filtrations}} \\
\hline
\textbf{eT} & Density of triangles \\
\textbf{eS} & Density of chordless squares \\
\textbf{eP} & Density of chordless pentagons \\
\textbf{e$\Sigma$} & Sum of cycle densities (triangles, squares, pentagons) \\
\textbf{eO} & Ollivier--Ricci curvature \\
\textbf{eF} & Augmented Forman--Ricci curvature \\
\textbf{eR} & Randi\'{c} connectivity index \\
\textbf{eH} & Harmonic index \\
\textbf{eA} & Repulsion-attraction rule \\
\textbf{eB} & Edge betweenness centrality \\
\hline
\multicolumn{2}{c}{\textbf{Node-based filtrations}} \\
\hline
\textbf{nD} & Degree \\
\textbf{nC} & Clustering coefficient \\
\textbf{nE} & Egonet persistence \\
\textbf{nG} & Graphlet-based \\
\hline
\multicolumn{2}{c}{\textbf{Metric-based filtrations}} \\
\hline
\textbf{mV} & Vietoris--Rips \\
\hline
\end{tabular}
\end{table}

\subsection{Densities of chordless cycles}

Given a graph $G$, let $e_{ij}$ be an edge between nodes $v_i$ and~$v_j$, with degrees $d_i>1$ and~$d_j>1$, respectively. The density of triangles corresponding to the edge~$e_{ij}$ is defined as the number $\#\triangle_{ij}$ of edge triangles in $G$ containing $e_{ij}$ divided by the maximum possible number of triangles in $G$ containing~$e_{ij}$, that~is,
\begin{equation}
\label{eq: density_triangles}
    C_t(e_{ij}) = \frac{\# \triangle_{ij}}{\min (d_i, d_j) - 1}.
\end{equation}

An edge cycle is called chordless if the only edges between its nodes are those forming the cycle.
Density of squares for $e_{ij}$ is defined by dividing the number of chordless squares in $G$ containing $e_{ij}$ by the maximum possible number of such squares given the 
degrees $d_i$ and $d_j$ and the existing triangles through~$e_{ij}$ in~$G$.
Density of pentagons for $e_{ij}$ is defined by counting the number of chordless pentagons containing $e_{ij}$ and normalizing by the maximum possible number of such pentagons, assuming known the degrees of $v_i$ and $v_j$ and the degrees of their respective neighbors.

\subsection{Experiments}

\vspace*{-0.2cm}

To gain insight into the behavior of chordless-cycle filtrations, we analyzed them as edge weights assigned to graphs within the same dataset. We observed that graphs belonging to different isomorphism classes exhibit distinct arrangements of edge weights. In contrast, some commonly used filtrations, such as Ollivier--Ricci curvature, produce nearly homogeneous edge weights on these datasets, failing to distinguish graphs already at this level. Consistent with this observation, we examined correlations between curvature-based measures and chordless-cycle densities in real-world networks, and found high correlations in the case of triangles (reaching a Pearson coefficient of $\rho=0.72$ for Ollivier--Ricci curvature in the wiki\_science network), and lower correlations for squares and pentagons 
(Figs.~S1 and S2 in the Supplementary Information).
Moreover, while most edge weightings separate graphs in simpler datasets, curvature- and centrality-based weights lose discriminative power as graph complexity increases. These observations suggest that filtrations based on chordless cycles encode structure sensitive to differences between graphs sharing similar regularities. Visual examples are provided in the Supplementary Information (Figs.\;S3--S5).

We analyzed the empirical performance of the proposed filtrations for distinguishing non-isomorphic graphs. We used datasets containing cubic graphs, strongly-regular graphs \cite{McKay}, and minimal Cayley graphs, as well as benchmark datasets for graph-learning tasks (BREC) \cite{Wang2023}. We provide detailed definitions of these datasets in Section~\ref{subsec:datasets} in the Appendix.

We also performed a property-prediction task on the ogbg-molhiv dataset \cite{ogbg}. For this, we followed the same approach as in \cite{Ballester2024}, including all the filtrations discussed in this paper, and we extended it to predict the following properties of graphs: average closeness centrality, average degree, average clustering coefficient, degree heterogeneity, average betweenness centrality, minimum degree, maximum degree, and average shortest path length, as well as those suggested in~\cite{Ballester2024}, namely maximum radius and diameter among the radii and diameters of the connected components of each graph, and their girth. From the persistence diagrams of the graphs, we computed persistence images \cite{adams2017persistence}, which were the input to a Random Forest regression model. We adopted the official OGB scaffold split, comprising 80\% training, 10\% validation, and 10\% test graphs. Models were trained on the training set and evaluated on the held-out test set. 
For the diameter, radius, and girth, we report accuracy after rounding predictions to the nearest integer, as in \cite{bause2025maximally}. For the remaining properties, we report the mean absolute error. 

To assess sensitivity of the topological filtrations under structural perturbations, we conducted systematic experiments on Barab\'asi--Albert (BA) \cite{barabasi1999emergence} and Erd\H{o}s--R\'enyi (ER) \cite{erdds1959random} random graph models, as well as on a Watts--Strogatz (WS) \cite{watts1998collective} model, which has a non-vanishing clustering coefficient.
We followed an approach similar to the one described in~\cite{tantardini2019comparing}. To facilitate a better comparison, we constructed graphs with $n=100$ nodes and the same average degree. For BA graphs, we varied the attachment parameter $m \in \{2, 5, 10\}$, resulting in average degrees of 4, 10, and 20. For ER graphs, we used connection probabilities $p \in \{0.04, 0.1, 0.2\}$, yielding similar average degrees (for $n=100$). For WS graphs, we generated small-world networks with a rewiring probability $\beta = 0.3$ and neighborhood sizes $\langle s \rangle \in \{4, 10, 20\}$, matching the corresponding average degrees of the ER and BA graphs. For each graph configuration across all models, we generated 50 independent instances and computed their persistence diagrams using 15 different filtrations. To simulate structural perturbations, we applied two modification schemes: edge random rewiring, which preserves node degrees but alters local structure, and edge random removal, which reduces connectivity. For each perturbation step, we recalculated the persistence diagram and measured the bottleneck distance to the original diagram in homology dimension $k = 0$. This process was repeated for 50 perturbation steps per graph, and results were averaged across 50 runs. 

\section{Results}
\subsection{Graph isomorphism problem}
\label{results}
We evaluate cubic, minimal Cayley, and strongly regular graphs, which are increasingly challenging benchmarks on which Weisfeiler--Leman tests are known to fail~\cite{morris2023weisfeiler}.
A~full description of the 15 tested filtrations is provided in Section\!~\ref{subsec:def_filtrations} in the Appendix. Throughout the text, we refer to them by their abbreviations stated in Table~\ref{tab:filtration_legend}.

Table~\ref{tab:ssr_k3} shows the success rate for distinguishing pairs of non-isomorphic strongly regular graphs when using different filtration weightings at clique expansion level $k=3$. The eS, eP, and e$\Sigma$ filtrations separate non-isomorphic graphs with a higher success rate than any other filtration weighting in our study. Strongly regular graphs form a particularly challenging dataset. It was found in~\cite{Ballester2024} that for $k=1$ most classical and curvature-based filtrations fail to distinguish non-isomorphic pairs.  In contrast, as shown in Table~\ref{tab:ssr_k3}, eP and e$\Sigma$ achieve high success rates on strongly regular graphs already at $k=1$, and perfect performance for $k=2$ and $k=3$. While the performance of eB filtration improves with higher clique expansions, eS, eP, and e$\Sigma$ consistently outperform all others across expansion levels, highlighting the importance of higher-order motif information.

\begin{table*}[htbp]
\renewcommand{\arraystretch}{1.3}
\caption{Success rate for distinguishing pairs of non-isomorphic strongly regular graphs when using different filtrations at clique expansion level of the graph with $k=3$. Bold indicates the best-performing algorithm. Filtrations are denoted by their abbreviations; see Table~\ref{tab:filtration_legend} for the complete legend.}
\label{tab:ssr_k3}
\centering
\scalebox{1.05}{
\begin{tabular}{l c c c c c c c c c c c c c c c c}
\toprule
Data & nD & eO & eF & mV & eT & eS & eP & e$\Sigma$ & eR & eH & eA & eB & nC & nE & nG\\
\midrule
\texttt{sr16622} & $\mathbf{1.00}$ & $\mathbf{1.00}$ & $\mathbf{1.00}$ & $\mathbf{1.00}$ & $\mathbf{1.00}$ & $\mathbf{1.00}$ & $\mathbf{1.00}$ & $\mathbf{1.00}$ & $\mathbf{1.00}$ & $\mathbf{1.00}$ & $\mathbf{1.00}$ & $\mathbf{1.00}$ & $\mathbf{1.00}$ & $\mathbf{1.00}$ & 
$\mathbf{1.00}$  \\
\texttt{sr251256} & $0.90$ & $0.90$ & $0.90$ & $0.90$ & $0.90$ & $\mathbf{1.00}$ & $\mathbf{1.00}$ & $\mathbf{1.00}$ & $0.90$ & $0.90$ & $0.90$ & $0.90$ & $0.90$ & $0.90$ & 
$\mathbf{1.00}$ \\
\texttt{sr261034} & $0.93$ & $0.93$ & $0.93$ & $0.93$ & $0.93$ & $0.96$ & $\mathbf{1.00}$ & $\mathbf{1.00}$ & $0.93$ & $0.93$ & $0.93$ & $0.93$ & $0.93$ & $0.93$ & 
0.96\\
\texttt{sr281264} & $\mathbf{1.00}$ & $\mathbf{1.00}$ & $\mathbf{1.00}$ & $\mathbf{1.00}$ & $\mathbf{1.00}$ & $\mathbf{1.00}$ & $\mathbf{1.00}$ & $\mathbf{1.00}$ & $\mathbf{1.00}$ & $\mathbf{1.00}$ & $\mathbf{1.00}$ & $\mathbf{1.00}$ & $\mathbf{1.00}$ & $\mathbf{1.00}$ & 
$\mathbf{1.00}$  \\
\texttt{sr291467} & $0.77$ & $0.77$ & $0.77$ & $0.77$ & $0.77$ & $\mathbf{1.00}$ & $\mathbf{1.00}$ & $\mathbf{1.00}$ & $0.77$ & $0.77$ & $0.77$ & $0.77$ & $0.77$ & $0.77$ & 
$\mathbf{1.00}$ \\
\texttt{sr351668} & $0.95$ & $0.95$ & $0.95$ & $0.95$ & $0.95$ & $\mathbf{1.00}$ & $\mathbf{1.00}$ & $\mathbf{1.00}$ & $0.95$ & $0.95$ & $0.95$ & $0.95$ & $0.95$ & $0.95$ & 
$\mathbf{1.00}$ \\
\texttt{sr351899} & $0.81$ & $0.81$ & $0.81$ & $0.81$ & $0.81$ & $\mathbf{1.00}$ & $\mathbf{1.00}$ & $\mathbf{1.00}$ & $0.81$ & $0.81$ & $0.81$ & $0.81$ & $0.81$ & $0.81$ & 
$\mathbf{1.00}$ \\
\texttt{sr361446} & $0.92$ & $0.92$ & $0.92$ & $0.92$ & $0.92$ & $0.98$  & $\mathbf{1.00}$ & $\mathbf{1.00}$ & $0.92$ & $0.92$ & $0.92$ & $0.92$ & $0.92$ & $0.92$ & 
0.98 \\
\texttt{sr401224} & $0.94$ & $0.94$ & $0.94$ & $0.94$ & $0.94$ & $0.98$ & $\mathbf{1.00}$ & $\mathbf{1.00}$ & $0.94$ & $0.94$ & $0.94$ & $0.94$ & $0.94$ & $0.94$ & 
0.99 \\
\bottomrule
\end{tabular}
}
\end{table*}

Tables~S7--S9 in the Supplementary Information consistently show a better performance of e$\Sigma$ across different expansion levels and datasets.
Our results on cubic graphs and minimal Cayley graphs are reported in the Supplementary Information as well (Tables~S3--S6).

For cubic graphs, nD fails at $k=1$ due to regularity, but improves at $k=2$ when higher-order cliques are included. At $k=1$, e$\Sigma$ and eB already achieve strong performance, outperforming mV and classical index-based filtrations, like eR, eH, eA, nC or nE. At $k=2$, eO also becomes effective, though cycle-based filtrations (eS, eP and e$\Sigma$) remain among the most reliable. These results underscore the importance of incorporating higher-order motifs for distinguishing highly regular graphs.

On minimal Cayley graphs, eB achieves perfect accuracy across expansion levels, reflecting sensitivity to global connectivity; and e$\Sigma$ also performs exceptionally well. Distance-based filtrations, like mV, and degree-based, such as nD, underperform, indicating that information based only on local or metric features is insufficient to capture Cayley graph structure.

In Table~\ref{tab:sota_brec_table}, we report the success rate for distinguishing pairs of instances of the BREC benchmark dataset~\cite{Wang2023} when using different filtrations with $k = 1$ and $k = 4$, and we also provide a comparison with state-of-the-art results. Bold indicates the best performing algorithm; $3$-WL refers to the Weisfeiler--Leman test \cite{morris2019weisfeiler}, $N_2$ corresponds to the $N_2$ algorithm \cite{papp2022theoretical}, and Ego to the egonet distance \cite{piccardi2023metrics}, considered among the best state-of-the-art techniques for the graph isomorphism problem (see~Section\!~\ref{sec:comparison_ego_dist} in the Appendix). We further evaluate all filtrations on the BREC dataset in Tables~S10--S13 in the Supplementary Information, including six graph families of varying structural complexity, and we report average performance across all categories. As in~\cite{Ballester2024}, due to combinatorial constraints, we did not calculate the Vietoris--Rips filtration for $k = 4$. We also report the success rate for distinguishing pairs of non-isomorphic cubic, minimal Cayley, and strongly regular graphs when using the ego-distance, in Table~S1 in the Supplementary Information.

On the BREC dataset, mV, nD, eR, and eH filtrations consistently perform worst across all $k$. Among curvature-based methods, eO outperforms eF but remains inferior to cycle-based filtrations like eS, eP, and e$\Sigma$. For $k=1$, e$\Sigma$ is the most effective filtration, followed by eS and eB. At higher expansion levels, eB becomes competitive, but cycle-density filtrations (eT, eS, eP, and e$\Sigma$) remain the most stable and consistently high-performing across all $k$.
Overall, these results show the high expressivity of persistent homology when an appropriate filtration is selected. 
Our e$\Sigma$ filtration not only achieves a high success rate in distinguishing pairs of graphs from the BREC dataset, but also performs consistently across different clique expansion levels $k$ of the graph. 
This makes our filtration a strong baseline for graph isomorphism and graph learning tasks.

Together, these results show that the choice of filtration largely determines the effectiveness of persistent homology for graph isomorphism. Across cubic, minimal Cayley, strongly regular, and BREC benchmarks, cycle-based filtrations perform best overall. Their advantage is strongest on difficult regular graph families, where filtrations based on degree, distance, and curvature often fail, especially at low clique levels. Although increasing \(k\) improves several methods, e$\Sigma$ remains the most accurate and stable across datasets, highlighting the value of higher-order motif information and making it a strong baseline for graph isomorphism and graph learning.

\begin{table*}[htbp]
\renewcommand{\arraystretch}{1.3}
\caption{Success rate for distinguishing pairs of instances of the BREC data set when using different filtrations with $k = 1$ and $k = 4$, and comparison with state-of-the-art (SOTA) results. Bold indicates the best-performing algorithm. $3$-WL refers to the Weisfeiler--Leman test \cite{morris2019weisfeiler}, $N_2$ corresponds to the $N_2$ algorithm \cite{papp2022theoretical}, and Ego to the egonet distance \cite{piccardi2023metrics}.
Filtrations are denoted by their abbreviations; see Table~\ref{tab:filtration_legend} for the complete legend.}
\label{tab:sota_brec_table}
\centering
\scalebox{1.05}{
\begin{tabular}{l ccc cccc cccc}
\toprule
Data & \multicolumn{3}{c}{SOTA} & \multicolumn{4}{c}{Filtration ($k = 1$)} & \multicolumn{4}{c}{Filtration ($k = 4$)} \\
\cmidrule(lr){2-4} \cmidrule(lr){5-8} \cmidrule(lr){9-12}
& 3-WL & $N_2$ & Ego & nD & eO & e$\Sigma$ & eB & nD & eO & e$\Sigma$ & eB \\
\midrule
Basic        & \textbf{1.00} & \textbf{1.00} & 0.76 & 0.03 & 0.93 & 0.98 & \textbf{1.00} & 0.83 & \textbf{1.00} & \textbf{1.00} & \textbf{1.00} \\
Regular      & 0.36 & 0.99 & 0.63 & 0.00 & 0.42 & \textbf{1.00} & 0.50 & 0.78 & 0.84  & \textbf{1.00} & 0.97 \\
Extension    & \textbf{1.00} & \textbf{1.00} & 0.73 & 0.07 & 0.76 &  0.95 & \textbf{1.00} & 0.29 & 0.92 &  0.98 & \textbf{1.00} \\
CFI          & 0.60 & 0.00 & \textbf{0.90} & 0.03 & 0.03 & 0.03 & 0.10 & 0.03 & 0.03 & 0.03 & 0.10 \\
\bottomrule
\end{tabular}
}
\end{table*}

\subsection{Property prediction}

\begin{figure}[htbp]
    \centering
\includegraphics[width=\columnwidth]{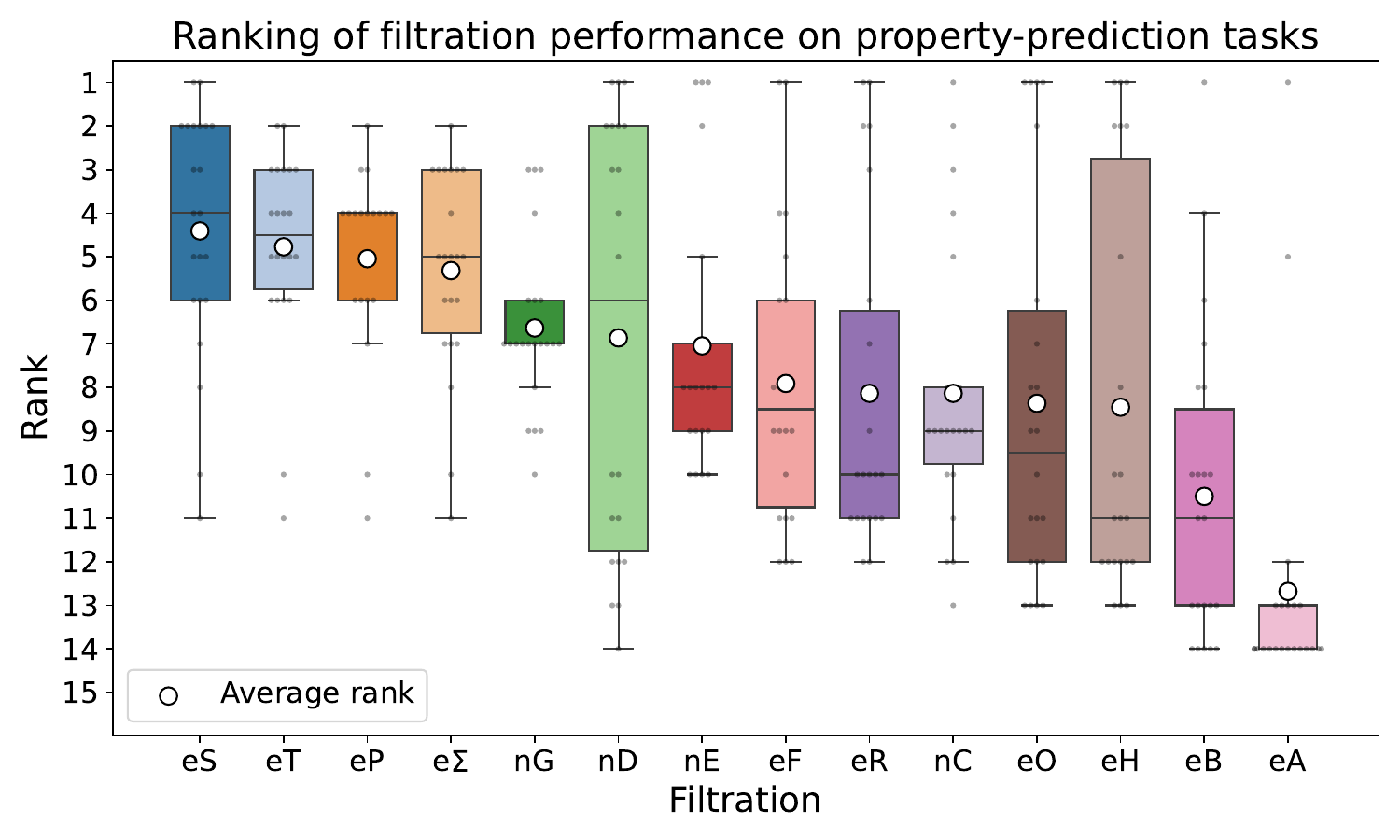}
    \caption{Distribution of ranks for each filtration method across graph property-prediction tasks. For each property and dimension (with $k=1$ and $k=2$), the filtrations are ranked by predictive performance. Ranks are computed for each prediction task, and then averaged across all tasks, so that tied methods share the same rank. Boxplots summarize the distribution of ranks across all tasks; gray points represent individual rank values, and white dots indicate the average rank for each filtration. Lower ranks correspond to better average predictive performance. Filtrations are denoted by their abbreviations; see Table~\ref{tab:filtration_legend} for the complete legend.}
    \label{fig:prop_pred_ranking}
\end{figure}

Figure\!~\ref{fig:prop_pred_ranking} summarizes the results of the property-prediction experiments, where the goal is to predict several structural properties of graphs: radius, diameter, girth, average closeness centrality, average degree, average clustering coefficient, degree heterogeneity, average betweenness centrality, minimum degree, maximum degree, and average shortest path length; using the ogbg-molhiv molecular graph dataset \cite{ogbg} under different filtrations. For each property and dimension, the filtrations are ranked by predictive performance. For radius, girth, and diameter, higher accuracy corresponds to better performance and therefore ranks are assigned in descending order; for all other properties, lower mean absolute error corresponds to better performance and ranks are assigned in ascending order. Ranks are computed independently for each prediction task using the ``min'' ranking method so that tied methods share the same rank, and then aggregated across all tasks. The boxplots summarize the distribution of ranks for each filtration across all prediction tasks, with lower ranks indicating better average predictive performance. Further details of the accuracy and mean average error results are provided in Tables~S17--S20 (Supplementary Information).

The boxplots for the property prediction task in Fig.\;\ref{fig:prop_pred_ranking} show that eS, eT, eP, and e$\Sigma$ achieve consistently low ranks, indicating the strongest and most stable predictive performance across properties.

Property-prediction tasks provide a more sensitive and interpretable evaluation of filtration expressivity, as they directly probe how well different topological summaries capture intrinsic graph structure. In the prediction tasks, the cycle-based filtrations (eS, eT, eP, and e$\Sigma$) consistently achieve the best rankings overall, with low median ranks and relatively small interquartile ranges, indicating stable performance across tasks. In contrast, eB and especially eA perform markedly worse. The filtrations nD, nE, eF, and eR exhibit intermediate behavior, but with noticeably wider interquartile ranges in some cases, suggesting that their performance is more task-dependent. These results support the use of cycle-density filtrations as descriptors for learning tasks, as they, overall, dominate the benchmark. 

This strong performance can be attributed to the fact that many of the predicted properties, such as diameter, girth, average shortest-path length, and average clustering coefficient, are intrinsically linked to the presence and distribution of cycles across multiple scales. Cycle-density filtrations explicitly encode connectivity patterns that are invisible to other filtrations. When combined with persistent homology, these filtrations retain information about the full distribution of topological features rather than collapsing them into single summary statistics, leading to improved and more stable predictive performance.

\subsection{Sensitivity tests}
Sensitivity analysis provides insight into how different filtrations respond to structural perturbations of a graph. Depending on the type of information they encode, some filtrations may remain stable under small perturbations, while others react sharply to changes in connectivity. Persistent homology offers a way to quantify these responses by measuring changes in persistence diagrams under different perturbations.

\begin{figure*}[htbp]
    \centering
\includegraphics[width=0.9\textwidth]{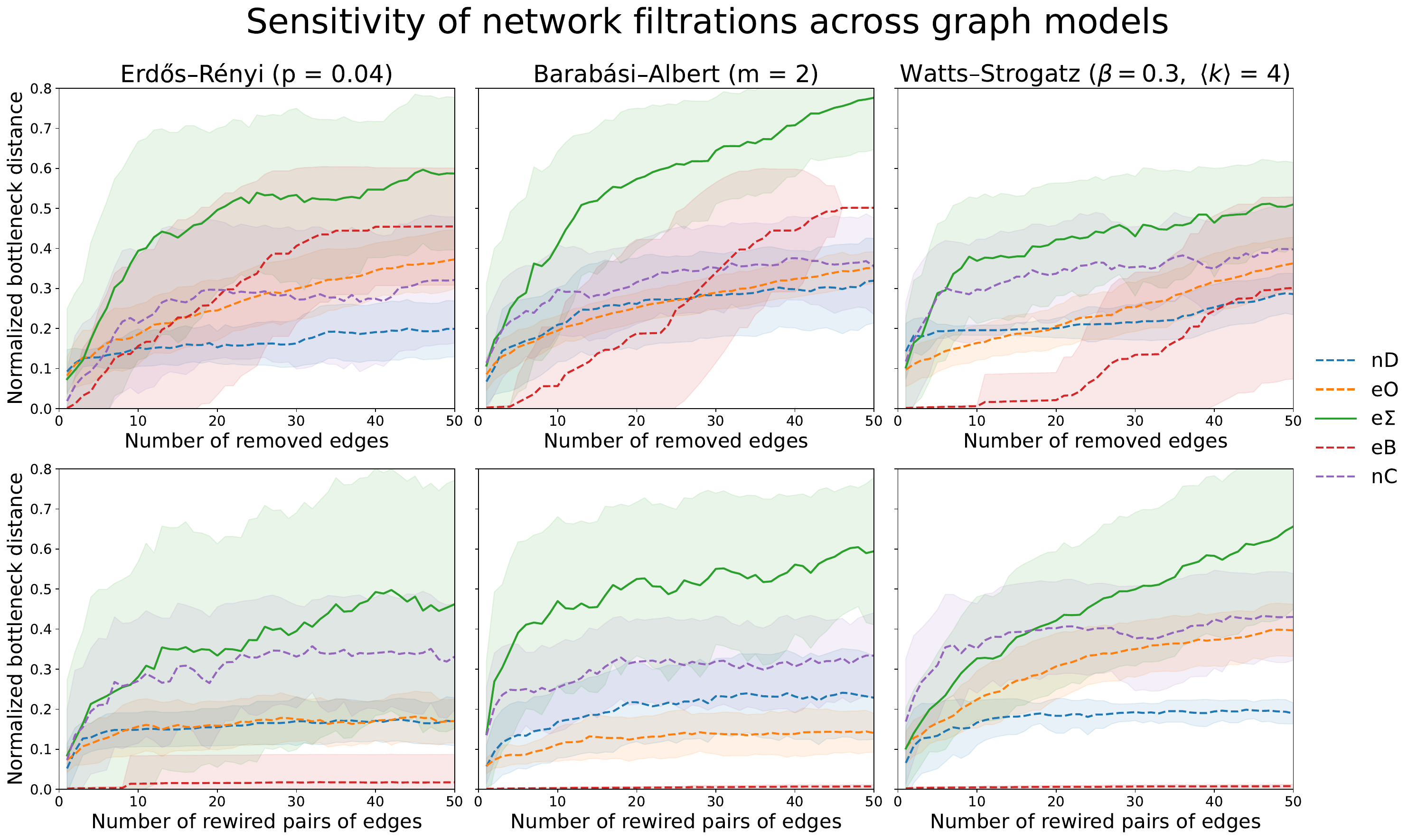}
\caption{Sensitivity of network filtrations across random graph models under edge removal (upper row) and edge rewiring (lower row) perturbations. Each panel shows the average normalized bottleneck distance (in homology dimension $k = 0$), computed over 50 runs of independently generated graphs with $n = 100$, as a function of the number of removed edges. Shaded regions indicate standard deviation. Results are presented for ER graphs with connection probability $p = 0.04$, BA graphs with attachment parameter $m = 2$, and WS small-world graphs with rewiring probability $\beta = 0.3$ and neighborhood sizes $\langle s \rangle = 4$. All configurations yield an approximate average degree of~4. The performance of the following filtrations, denoted by their abbreviations, is reported: nD, eO, e$\Sigma$, eB, and~nC; see Table~\ref{tab:filtration_legend} for the complete legend.}
\label{fig:sensitivity}
\end{figure*}

A comparison of the sensitivity of the different filtrations to edge removal and edge rewiring perturbations across different random graph models is shown in Fig.\;\ref{fig:sensitivity}. In particular, we select the best-performing and most significant filtrations that represent different types of measures (centrality, edge-based, degree, cycle-density, etc.). We compute the normalized bottleneck distance by dividing the filtration values in each persistence diagram by the maximum filtration value, following the definition in \cite{may2024normalized}. This normalization ensures that bottleneck distances are comparable across filtrations with different value ranges, allowing the observed differences in sensitivity to be attributed to structural expressivity rather than scale effects. Among the tested filtrations, e$\Sigma$ followed by nC exhibits the highest sensitivity to edge removal and edge rewiring. Both filtrations are particularly responsive to changes in local connectivity patterns. 
As one would expect, under degree-preserving randomization (edge rewiring), nD shows smaller changes in the bottleneck distance, whereas e$\Sigma$ and nC, based on cycles and clustering, respectively, detect topological modifications earlier. 

Beyond overall sensitivity, the temporal profile of the response is also informative. Here e$\Sigma$ and nC exhibit a sharp increase in bottleneck distance even for a small number of perturbations, indicating strong sensitivity to early structural changes. In contrast, nD often presents a gradual response. This suggests that cycle-density filtrations can act as early detectors of topological change. 

In the Supplementary Information (Figs.\;S6 and S7), we extend the sensitivity analysis to different parameters of the same random graph models, and we show the results for all the filtrations used in this paper.
The results include Erd\H{o}s–R\'enyi (ER) graphs with connection probabilities $p \in \{0.04, 0.1, 0.2\}$, Barab\'asi–Albert (BA) graphs with attachment parameters $m \in \{2, 5, 10\}$, and Watts–Strogatz (WS) small-world graphs with rewiring probability $\beta = 0.3$ and with neighborhood sizes $\langle s \rangle \in \{4, 10, 20\}$, such that all configurations yield comparable average degrees across models. We observe that the normalized bottleneck distance generally increases with the number of removed or rewired edges, indicating reduced topological stability under progressive loss or change of connectivity. This effect is most pronounced in sparse graphs (e.g., ER with $p = 0.04$ and BA with $m = 2$), while denser networks show lower sensitivity. From a percolation perspective, this behavior indicates that cycle-based filtrations are sensitive to changes in path redundancy, not just connectivity, and detect the progressive loss of alternative routes well before the network reaches a macroscopic connectivity breakdown.

\section{Discussion}
In this article, we use a new class of filtrations for persistent homology, based on densities of edge triangles and chordless edge squares and edge pentagons, which were introduced in~\cite{almagro2022detecting} for dimensionality detection purposes. Our results show that, even at clique level $k=1$, these filtrations are highly expressive in distinguishing non-isomorphic graphs. Moreover, even if performance improves for higher values of the clique simplicial dimension~$k$, being effective for $k=1$ helps reducing time complexity.

The term graph expressivity, as in \cite{Ballester2024}, refers to two related but distinct concepts: the ability of a method to distinguish non-isomorphic graphs, and the ability to capture certain graph properties.
Beyond graph isomorphism, the property-prediction experiments provide additional insight into the amount of structural information captured by cycle-based filtrations. 
While graph isomorphism is highly sensitive to structural perturbations, since removing a single edge suffices to alter the isomorphism type, property prediction is feasible for real-world networks and more directly relevant to applications. 

The strong and consistent performance of eS, eT, eP, and e$\Sigma$ suggests that these filtrations encode mesoscale organization, bridging purely local descriptors, like degree, and global summaries, such as spectra. By weighting edges according to their participation in cycles, motif-based filtrations offer a representation that aligns with the structural mechanisms of global graph properties. 

The fact that cycle-based filtrations are very sensitive to small edge perturbations also explains their ability to detect small structural changes, and complements the isomorphism and property-prediction results. Filtrations that respond early to edge modifications are better suited to capture subtle but meaningful differences in connectivity patterns, thus distinguishing highly similar but non-isomorphic graphs. This explains the strong performance of the chordless-cycle filtration in challenging settings.

Our motif-based filtrations consistently achieve high discrimination across highly regular graph families such as cubic and minimal Cayley graphs, whereas degree-based filtrations fail. Strongly regular graphs are notoriously difficult for isomorphism detection. Our e$\Sigma$ filtration is one of the few methods achieving perfect distinction at higher expansion levels ($k = 2,\ 3$), and very high success even at $k = 1$, far outperforming eB and the ego-distance method. Although the eB filtration was highly competitive on simpler datasets, it largely fails here, suggesting that shortest-path information alone is insufficient for such symmetric structures. And computing it for large graphs is computationally expensive. Moreover, nE
underperforms the e$\Sigma$ filtration across all $k$. This suggests that capturing local cyclic motifs, such as chordless cycles, provides discriminative power when other metrics collapse due to global symmetry. Similarly, the e$\Sigma$ filtration consistently outperforms all other methods across various subsets and across different values of~$k$ for the BREC dataset, and it remains the top-performing method at low clique levels. The eB filtration becomes competitive for $k \geq 2$, but does not outperform e$\Sigma$ overall.

Filtrations based on simple local features, such as nD, eR, or eH, consistently underperform, particularly on highly symmetric or regular graphs. These measures depend mainly on the degrees of nodes and ignore deeper connectivity or motif structure, being ineffective when graphs are globally different but locally similar. In contrast, motif-based filtrations, and, in particular, eS, eP, and e$\Sigma$, capture local subgraph patterns of increasing complexity and are still discriminative in highly symmetric settings. The nG filtration also achieves excellent performance as it encodes node information in small induced subgraphs. However, it requires counting orbits, which is computationally expensive, especially compared to our e$\Sigma$ approach. Centrality-based methods like eB are effective on graphs where shortest-path structure diverges but fail in highly regular graphs. 

Our filtrations are computationally scalable. The triangle, square, and pentagon filtrations require $O(|E|\Delta^2)$, $O(|E|\Delta^3)$, and $O(|E|\Delta^4)$ time, respectively, where $|E|$ is the number of edges and $\Delta$ the maximum node degree. For sparse real-world networks, where $|E| = O(|N|)$ and the degrees are typically bounded, these complexities become linear in the number of nodes $|N|$ of the graph. Persistent homology calculations on simplicial complexes, such as clique complexes of graphs, can be performed in at most $O(n^3)$ time, where $n$ is the number of simplices~\cite{EdelsbrunnerHarer2010}. If only the underlying graphs are considered (i.e., at clique dimension $k=1$), then the bound can be reduced to near-linear time on the number of filtration steps~\cite{DeyHouParsa2023}.

A natural question is whether edge weightings alone are enough to distinguish non-isomorphic pairs of graphs, or whether persistent homology provides additional value. Empirical evidence supports that while weightings encode local information, persistent homology captures how such information is organized across scales. In~\cite{ferra2025}, the authors compared the performance of a neural network classifier using chordless cycle weightings and the same filtration combined with persistent homology, and found that the performance increases when both types of information are combined. 

Triangle density weightings are closely correlated with curvatures, as displayed in Figs.\;S1 and S2 in the Supplementary Information. In \cite{Jost-Liu_2014}, the Ollivier--Ricci curvature at an edge was estimated in terms of the degrees of its endpoints and the number of triangles containing the edge, thus obtaining lower bounds in terms of the clustering coefficients of the endpoints. Likewise, the augmented Forman--Ricci curvature at an edge incorporates the number of triangles to which the given edge belongs as an additive correction to the classical Forman--Ricci formula \cite{Sreejith_etal_2016,Samal2018,Fesser_2024}. Hence, chordless-cycle densities capture geometric information that extends the role of triangles in curvature to polygons of larger size.

While motif-based filtrations exhibit strong and consistent performance, their advantages are most pronounced when structural distinctions arise from higher-order connectivity rather than from degree-based characteristics. In graphs dominated by long-range correlations or with few cycles, distance-based methods may remain competitive. Moreover, although our approach is computationally efficient compared to graphlet-based methods, counting larger chordless cycles may become costly for very dense graphs, suggesting a trade-off between expressivity and scalability.

Our findings align with prior results in \cite{Ballester2024} and show that our method achieves and exceeds state-of-the-art performance while significantly reducing complexity. Hence, our results position the e$\Sigma$ filtration as a highly expressive, interpretable, and efficient approach for graph comparison tasks, well-suited for practical applications in domains such as neuroscience, chemistry, and network science. Moreover, our filtration could be adapted to measure dissimilarities between real-world networks.

\vfill

\begin{appendices}
\section*{APPENDIX}
\label{sec:appendix}

\subsection{Filtrations}
\label{subsec:def_filtrations}

A node filtration on a graph $G=(N,E)$ is a function $w_N\colon N\to\mathbb{R}$ and an edge filtration is a function $w_E\colon E\to\mathbb{R}$. A~node filtration $w_N$ and an edge filtration $w_E$ are compatible if $w_N(v)\le w_E(e)$ whenever $v$ is incident with~$e$.
 
The node and edge filtrations used in this study are the following:
\begin{enumerate}
    \item Degree filtration (denoted by \textbf{nD}): $v \mapsto \deg(v)$.  
\item Filtration based on Ollivier--Ricci curvature~\cite{ollivier2007ricci} (denoted by \textbf{eO}). We assign weight $-1$ to nodes and
\[
    \kappa(u, v) = 1 - W_1(\mu^\alpha_u, \mu^\alpha_v)
\]
to edges, where $ W_1 $ denotes the first Wasserstein distance \cite{mileyko2011probability} and $ \mu^\alpha_u$, $\mu^\alpha_v $ are probability measures based on a lazy random walk in the graph:
\[
\mu^\alpha_i(j) = 
\begin{cases} 
\alpha & \text{if } j = i, \\
(1 - \alpha)/\deg(i) & \text{if } i \sim j, \\
0 & \text{otherwise,}
\end{cases}
\]
where $\alpha$ is as a smoothing parameter, which we set to $\alpha = 0$ for our experiments, consistently with~\cite{Ballester2024}. 

\item Filtration based on the augmented Forman--Ricci curvature~\cite{Samal2018} (denoted by \textbf{eF}). We set $v \mapsto -1$ and $(u, v) \mapsto 4 - \deg(u) - \deg(v) + 3 |N(u) \cap N(v)|$, where $N(u)$ and $N(v)$ denote the neighborhoods of $u$ and~$v$, respectively.
\item Vietoris--Rips filtration (denoted by \textbf{mV}),
where distance is given by the shortest-path length between nodes~\cite{Adams2022}.
This filtration is neither determined by node weights nor by edge weights, but the parameter value associated with each simplex $\sigma$ in the Vietoris--Rips complex is the distance value at which $\sigma$ appears.
    \item Filtration assigning to each edge the density of triangles (denoted by \textbf{eT}), computed by Eq.\;\ref{eq: density_triangles}.
    \item Filtration assigning to each edge the density of chordless squares (denoted by \textbf{eS}). 
    \item Filtration assigning to each edge the density of chordless pentagons (denoted by \textbf{eP}). 
    \item Filtration assigning to each edge the sum of densities of triangles, squares, and pentagons (denoted by \textbf{e\boldmath$\Sigma$}).
    \item Filtration assigning to each edge the Randi\'{c} connectivity index \cite{randic1975characterization} (denoted by \textbf{eR}), which is defined as shown in Eq.\;\ref{eq:randic_index},     where $d_i$ and $d_j$ are the degrees of nodes $i$ and~$j$:
\begin{equation}\label{eq:randic_index}
        R(e_{ij})= 
        \frac{1}{\sqrt{d_i d_j}},
\end{equation}
   \item Filtration that assigns to each edge the weights used to compute the harmonic index \cite{fajtlowicz1987conjectures} (denoted by \textbf{eH}). The harmonic index of a graph $G$ is the sum of weights 
    $2/(d_u+d_v)$
    of all edges $(u,v)$ of $G$, where $d_u$ denotes the degree of a node $u$ in $G$.
    \item Filtration that assigns to each edge a variant of the  
    repulsion-attraction rule \cite{muscoloni2017machine} (denoted by \textbf{eA}):  \begin{equation}\label{eq:repulsion_attraction rule}
         \text{RA}
         (e_{ij}) = \frac{d_i + d_j + d_id_j}{1+|N(i)\cap N(j)|},
    \end{equation}
    where $d_i$ and $d_j$ denote the degrees of nodes $i$ and $j$, and 
    $|N(i)\cap N(j)|$ is the number of common neighbors between them.
    \item Filtration assigning to each edge its betweenness centrality \cite{muscoloni2017machine} (denoted by \textbf{eB}). The betweenness centrality of an edge $e_{ij}$ is computed as 
\begin{equation}\label{eq:betweenness_centrality}
c_B(e_{ij}) = \sum_{s,t \in N} \frac{\sigma(s,t \mid e_{ij})}{\sigma(s,t)}, 
    \end{equation}
    where $N$ is the set of nodes, $\sigma(s,t)$ is the number of shortest paths between $s$ and $t$, and $\sigma(s,t \mid e_{ij})$ is the number of those paths passing through $e_{ij}$.
    \item Clustering coefficient  filtration (denoted by \textbf{nC}). 
    If $T(v)$ is the number of triangles containing a node $v$ and $d_v$ is the degree of~$v$, then the clustering coefficient \cite{watts1998collective} is  
    $C(v) = 2T(v)/(d_v(d_v-1))$.
    \item The egonet $E_i$ of node~$i$ is the subgraph on nodes $\{i\}\cup N(i)$, where $N(i)$ is the set of neighbors of~$i$ \cite{piccardi2023metrics}.
    The egonet filtration (denoted by \textbf{nE}) assigns to each node the egonet persistence, defined for undirected and unweighted networks by \begin{equation}\label{eq:egonet_p}
        p(i) = \frac{\sum_{j \in E_i} m_j^{\mathrm{int}}}{\sum_{j \in E_i} (m_j^{\mathrm{int}} + m_j^{\mathrm{ext}})},
    \end{equation}
    where $m_j^{\mathrm{int}}$ and $m_j^{\mathrm{ext}}$ denote the internal and external degrees of $j$ in $E_i$ \cite{Piccardi2011}. 
    We set $p(i) = 0$ for isolated nodes.
    \item Graphlets \cite{sarajlic2016graphlet} are small, connected, non-isomor\-phic subgraphs. An orbit is a distinct node position within a graphlet, defined up to automorphism. For each node $v$, we compute its Graphlet Degree Vector (GDV), whose entries count how many times $v$ participates in each orbit. 
    The graphlet-based filtration (denoted by \textbf{nG}) assigns to $v$ the following score, where $i$ indexes the 11 graphlet orbits:
\[
\textstyle
\mathrm{graphlet\_score}(v) =
\sum_{i=0}^{10} 
\,\mathrm{o\_count}(v,i),
\]

\end{enumerate}

Filtrations nC and nE are included as comparison measures in our pipeline. We also evaluate Piccardi's original method~\cite{piccardi2023metrics} separately (Section\!~\ref{sec:comparison_ego_dist} and Table~\ref{tab:sota_brec_table}). 

\subsection{Datasets}\label{subsec:datasets}
Cubic graphs are graphs where each node has degree~$3$; such graphs are also called $3$-regular graphs. Some sources for connected cubic graphs are \cite{McKay}, \cite{Brinkmann2013}.
Minimal Cayley graphs \cite{Knauer2025} are Cayley graphs from finite groups $G$ equipped with a minimal set $S$ of generators, where minimality means that no proper subset of $S$ generates~$G$. The Cayley graph of $G$ relative to a generating set $S$ has the elements of $G$ as nodes and there is an edge from $x$ to $y$ if and only if $x^{-1}y\in S$.
A~strongly regular graph
$\text{SRG}(n,k,\lambda,\mu)$
has $n$ nodes, each of which has degree~$k$, and any two adjacent nodes have $\lambda$ common neighbors, while any two non-adjacent nodes have $\mu$ common neighbors \cite{McKay,Coolsaet2023House}. 

BREC is a benchmark dataset designed to evaluate graph distinguishability, consisting of 400 pairs of non-isomorphic graphs \cite{Wang2023}. These fall into four main categories: Basic (1-WL indistinguishable and non-regular); Regular (all nodes have the same degree); Extension (based on theoretical GNN extensions, whose distinction requires more expressive methods than standard GNNs \cite{papp2022theoretical}); and CFI (constructed with the Cai--Fürer--Immerman method, distinguishable by $k$-WL but not by $(k-1)$-WL  \cite{cai1992optimal}).
Following \cite{Ballester2024}, we subdivide some of these into additional families, including graphs satisfying the 4-vertex condition (with highly regular local 4\nobreakdash-node structure) and distance-regular graphs (with identical distance structure from every node).

\subsection{Experimental setup}\label{subsec:experimental_setup}
As in \cite{Ballester2024}, after selecting a filtration (excluding the Vietoris--Rips filtration), we expand the graph by adding cliques up to any chosen dimension $k+1$. The filtration value assigned to a clique is defined recursively as the maximum filtration value among its proper subcliques.
We then compute persistent homology up to dimension~$k$. For $k = 1$, no higher cliques are added, and we analyze the graph based solely on its connected components and cycles. 

The bottleneck distance is a common metric to compute similarity between two persistence diagrams. It is the shortest distance $d$ for which there exists a perfect matching between the points of the two persistence diagrams, completed with all the points on the diagonal in order to ignore cardinality mismatches, such that any couple of matched points are at distance at most~$d$, where the distance between points is the supremum norm in the plane
\cite{gudhi:FilteredComplexes}. 
Hence,
\begin{equation}\label{eq:bottleneck_distance}
    d_B(D, D') = \inf_{\eta : D \to D'} \sup_{x \in D} \|x - \eta(x)\|_\infty, 
\end{equation}
where $\eta$ ranges over all matchings between the two persistence diagrams $D$ and $D'$. Using our computations, we obtain several persistence diagrams for each graph. We compare them by pairs using the bottleneck distance, as defined in Eq.\;\ref{eq:bottleneck_distance}. We consider two graphs to be different when the distance between their persistence diagrams is $> 10^{-8}$, that is, above machine precision.

\subsection{Comparison with ego-distance}\label{sec:comparison_ego_dist}
We compare our approach with the egonet-based distance proposed in~\cite{piccardi2023metrics}. 
For each node~$i$, we computed three local descriptors: normalized degree $nd(i)$, local clustering coefficient $c_i$, and egonet persistence $p_i$ (Eq.\;\ref{eq:egonet_p}).
Following~\cite{piccardi2023metrics}, we discretize $[0,1]$ into $r$ bins and construct the joint 3D histogram $P_{d,c,p}$ of $(nd, c, p)$ values over all nodes. Let $Q_{d,c,p}$ denote the associated cumulative distribution.
The ego-distance between graphs $G'$ and $G''$ is defined as the $\ell_2$ norm of the difference between their cumulative distributions \cite{piccardi2023metrics}:
\begin{equation}
\label{eq:egonet_based_distance}
 \textstyle\left ( \sum_{h,k,n= 1}^r (Q'_{d,c,p}(h,k,n) - Q''_{d,c,p}(h,k,n))^2 \right )^\frac{1}{2}.
\end{equation}

We use a cap value $T = 0.5$ and a discretization step $\Delta = 0.01$, as recommended in~\cite{piccardi2023metrics}. 
We consider two graphs to be distinct
if their ego-distance exceeds $10^{-8}$, hence above machine precision. We report these results in Table~\ref{tab:sota_brec_table}.

\end{appendices}

\section*{ACKNOWLEDGMENTS}
M.V.M. acknowledges funding from the Luddy School of Informatics, Computing, and Engineering (Indiana University Bloomington) to attend NetSci 2025.
R.J. acknowledges support from the fellowship FI-SDUR funded by Generalitat de Catalunya. 
M.A.S. acknowledges support from grant  TED2021-129791B-I00 funded by MCIN/AEI/10.13039/501100011033 and by \emph{European Union NextGenerationEU/PRTR}, and grant PID2022-137505NB-C22 funded by MCIN/AEI/ 10.13039/501100011033 and by ERDF's \emph{A~way of making Europe}.
A.F.M., R.B., and C.C. were partially funded through MCIN/AEI
grants PID2020-117971GB-22 and PID2022-136436NB-I00, as well as AGAUR grant 2021 SGR 00697. A.F.M. was also supported by MCIN/AEI under grant PRE2020-094372, and R.B. was supported through the FPU contract PU21/00968.

\section*{DATA AND CODE AVAILABILITY}
The open-source code, along with the code to reproduce the figures, will be available on GitHub upon publication. The expressivity experiments using known graphs (cubic, minimal Cayley, and strongly regular) use data from the database in \cite{Coolsaet2023House}, which does not specify any licensing requirements but requires a citation. The BREC data set is distributed under an MIT license. We used part of the code from the authors of~\cite{Ballester2024}, which is available under a 3-Clause BSD License. The MATLAB code for the function EgoDist, which computes ego-distances, and its Python adaptation are available at~\url{https://piccardi.faculty.polimi.it/highlights.html}.

\section*{COMPETING INTERESTS STATEMENT}
The authors declare no competing interests.	

\vfill

\bibliography{references}

\onecolumngrid

\clearpage
\pagestyle{empty}



\section*{Supplementary material}

# Supplementary Information for "Motif-based filtrations for persistent homology: A framework for graph isomorphism and property prediction"

Meritxell Vila-Miñana,[1] Robert Jankowski,[2,3,4] Aina Ferrà Marcús,[5]
Rubén Ballester,[5] M. Ángeles Serrano,[3,4,6] and Carles Casacuberta[5]

[1] Center for Complex Networks and Systems Research, Luddy School of Informatics,
Computing, and Engineering, Indiana University, Bloomington, IN, USA
[2] Faculty of Electrical Engineering, Mathematics and Computer Science,
Delft University of Technology, 2628 CD, Delft, Netherlands
[3] Departament de Física de la Matèria Condensada,
Universitat de Barcelona, Martí i Franquès 1, E-08028 Barcelona, Spain
[4] Institute of Complex Systems (UBICS), Universitat de Barcelona, Barcelona, Spain
[5] Departament de Matemàtiques i Informàtica, Universitat de Barcelona,
Gran Via de les Corts Catalanes 585, E-08007 Barcelona, Spain
[6] ICREA, Passeig Lluís Companys 23, E-08010 Barcelona, Spain

**EXTENDED RESULTS**

Table S1 shows the success rate for distinguishing pairs of non-isomorphic graphs from the following data bases: cubic graphs, minimal Cayley graphs, and strongly regular graphs, and pairs of instances of the BREC data set using the ego-distance.

TABLE S1: Success rate for distinguishing pairs of non-isomorphic *cubic*, *minimal Cayley*, and *strongly regular* graphs, and pairs of instances of the BREC data set when using the ego-distance.

| Data  | Success Rate | Data  | Success Rate | Data     | Success Rate | Data | Success Rate |
|-------|--------------|-------|--------------|----------|--------------|------|--------------|
| cub06 | 0.00         | cay12 | 0.71         | sr16622  | 0.93         | B    | 0.76         |
| cub08 | 0.00         | cay16 | 0.73         | sr251256 | 0.93         | R    | 0.63         |
| cub10 | 0.00         | cay20 | 0.67         | sr261034 | 1.00         | E    | 0.73         |
| cub12 | 0.00         | cay24 | 0.75         | sr281264 | 0.93         | C    | 0.90         |
| cub14 | 0.00         | cay32 | 0.77         | sr291467 | 0.94         | 4    | 0.97         |
|       |              | cay36 | 0.77         | sr351668 | 1.00         | D    | 0.99         |
|       |              | cay60 | 0.75         | sr351899 | 0.99         |      |              |
|       |              | cay63 | 0.75         | sr361446 | 1.00         |      |              |
|       |              |       |              | sr401224 | 1.00         |      |              |

In Figures S1 and S2, we represent the relationship between the edge weights based on the Ollivier–Ricci and Forman–Ricci curvatures and density of chordless cycles of triangles, squares and pentagons, for some of the real networks described in Table S2. One can observe that the correlation between curvature measures and density of squares and pentagons is low. Hence, we can infer that they are capturing complementary information to each other, which makes chordless-cycle densities an appealing choice for a TDA filtration. Additionally, we note that the edge weighting based on the density of triangles is most similar to the Ollivier–Ricci curvature.

In Figures S3, S4, and S5, we show four different graphs from the cubic, minimal Cayley, and strongly regular datasets, respectively. Each row corresponds to a different measure of edge weights. The lighter the color, the larger the values are. To distinguish four non-isomorphic graphs, we note that their weight distributions differ. Qualitatively, we observe that for the cubic dataset, all filtrations reveal differences among the four graphs. In the case of minimal Cayley, Ollivier–Ricci curvature does not distinguish between graphs 2 and 4. Similarly, triangle weighting provides no differences between graphs, and pentagons do not distinguish between graphs 1 and 4. Instead, betweenness centrality and the cycle densities allow us to distinguish all of them. For strongly regular graphs, while Ollivier–Ricci curvature, triangles, and betweenness centrality weightings provide no differences between graphs, the squares, pentagons, and the combination of density of chordless cycles allow us to distinguish all of them.

Tables S3, S5 and S7 contain numerical results corresponding to this discussion. These show the success rates for distinguishing pairs of non-isomorphic cubic, minimal Cayley, and strongly regular graphs, respectively, when using

2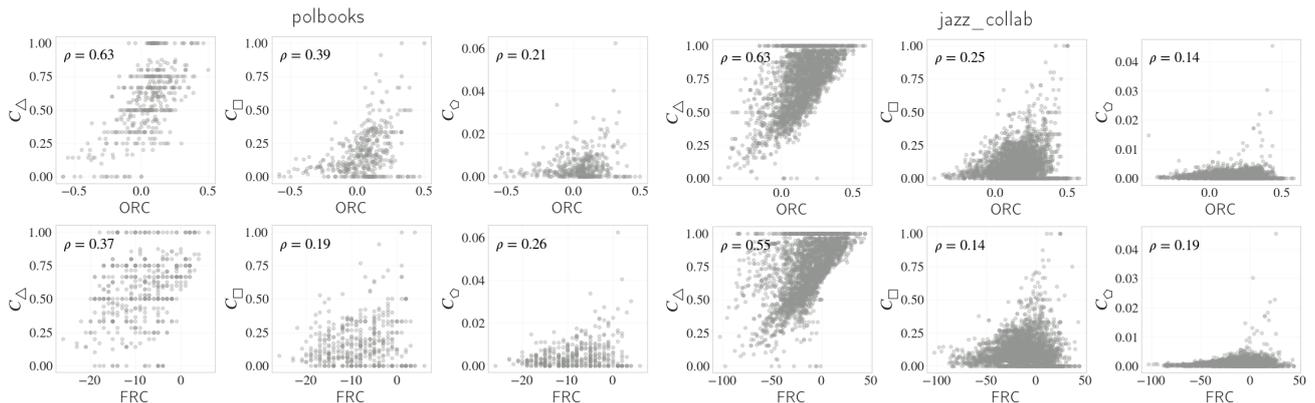

FIG. S1: Relationship between edge weights based on Ollivier–Ricci and Forman–Ricci curvatures and density of chordless cycles of triangles, squares and pentagons, for the real networks *polbooks* and *jazz_collab*.

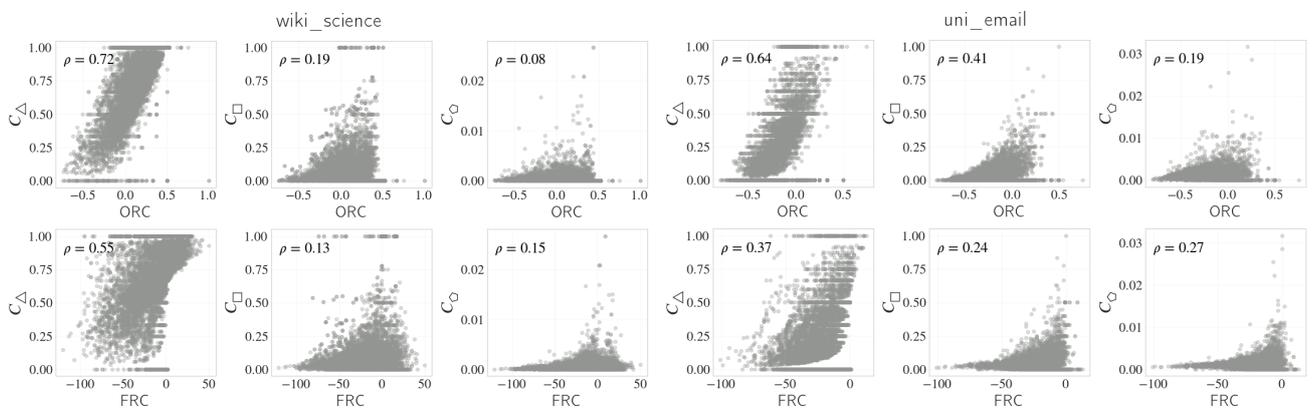

FIG. S2: Relationship between edge weights based on Ollivier–Ricci and Forman–Ricci curvatures and density of chordless cycles of triangles, squares and pentagons, for the real networks *wiki_science* and *uni_email*.

different filtrations at the clique-expansion level $k = 1$. Tables S3, S4 contain the success rate for distinguishing pairs of non-isomorphic cubic graphs when using different filtrations at expansion levels $k = 1$ and $k = 2$. Tables S5, S6 contain the success rate for distinguishing pairs of non-isomorphic minimal Cayley graphs when using different filtrations at expansion levels $k = 1$ and $k = 2$. Tables S7 and S8 contain the success rate for distinguishing pairs of non-isomorphic strongly regular graphs when using different filtrations at expansion levels $k = 1$ and $k = 2$. Tables S9, S10, S11, S12 contain the success rate for distinguishing pairs of instances of the BREC dataset when using different filtrations at expansion levels $k = 1$, $k = 2$, $k = 3$, and $k = 4$. This dataset [1] contains 400 graphs from six categories: Basic (60), Regular (100), Extension (100), CFI (100), 4-vertex Condition (20), and Distance-Regular (20).

Table S13 provides a summary of the statistics for all graph properties in the ogbg-molhiv data set [2], and Tables S14, S15, S16 and S17, respectively for dimensions $k = 1$ and $k = 2$, contain the results of the property-prediction experiments, consisting in predicting the following properties of graphs: radius, diameter, girth, average closeness centrality, average degree, average clustering coefficient, degree heterogeneity, average betweenness centrality, minimum degree, maximum degree, and average shortest path length, using the ogbg-molhiv graph molecular data set [2] and the different filtrations studied in this paper.

Finally, Figures S6 and S7 contain a comparison of the sensitivity of different filtrations to edge removal and edge rewiring perturbations across different random graph models and different parameters.

TABLE S2: Description of four real-world networks and their properties, including the number of nodes and edges.

| Name | Description | Nodes ($N$) | Edges ($E$) |
|---|---|---|---|
| *polbooks* | Political books network. A network of books about U.S. politics published close to the 2004 U.S. presidential election, and sold by Amazon.com. Edges between books represent frequent copurchasing of those books by the same buyers. | 105 | 441 |
| *jazz_collab* | Jazz collaboration network. The network of collaborations among jazz musicians, and among jazz bands, extracted from The Red Hot Jazz Archive digital database, covering bands that performed between 1912 and 1940. | 198 | 2742 |
| *wiki_science* | Wikipedia Map of Science. A network of scientific fields, extracted from the English Wikipedia in early 2020. Nodes are Wikipedia pages representing natural, formal, social, and applied sciences, and two nodes are linked if the cosine similarity of the page content is above a threshold. | 687 | 6523 |
| *uni_email* | Email network. A network representing the exchange of emails among members of the Rovira i Virgili University in Spain in 2003. | 1133 | 10903 |

TABLE S3: Success rate for distinguishing pairs of non-isomorphic *cubic* graphs when using different filtrations at expansion level of the graph $k = 1$. Bold indicates the best-performing algorithm. Filtrations are denoted by their abbreviations; see Table I in the main text for the complete legend.

| Data | nD | eO | eF | mV | eT | eS | eP | e$\Sigma$ | eR | eH | eA | eB | nC | nE | nG |
|---|---|---|---|---|---|---|---|---|---|---|---|---|---|---|---|
| cub06 | 0.00 | **1.00** | **1.00** | 0.00 | **1.00** | **1.00** | 0.00 | **1.00** | 0.00 | 0.00 | **1.00** | **1.00** | **1.00** | **1.00** | **1.00** |
| cub08 | 0.00 | 0.90 | 0.70 | 0.00 | 0.70 | **1.00** | 0.90 | **1.00** | 0.00 | 0.00 | 0.90 | **1.00** | 0.40 | 0.40 | **1.00** |
| cub10 | 0.00 | 0.98 | 0.66 | 0.00 | 0.66 | 0.89 | 0.91 | 0.99 | 0.00 | 0.00 | 0.85 | **1.00** | 0.11 | 0.11 | 0.98 |
| cub12 | 0.00 | 0.99 | 0.64 | 0.00 | 0.62 | 0.88 | 0.92 | **1.00** | 0.00 | 0.84 | **1.00** | 0.14 | 0.14 | 0.14 | 0.98 |
| cub14 | 0.00 | 0.99 | 0.62 | 0.00 | 0.61 | 0.85 | 0.94 | **1.00** | 0.00 | 0.83 | **1.00** | 0.12 | 0.12 | 0.12 | 0.98 |

TABLE S4: Success rate for distinguishing pairs of non-isomorphic *cubic* graphs when using different filtrations at expansion level of the graph $k = 2$. Bold indicates the best-performing algorithm. Filtrations are denoted by their abbreviations; see Table I for the complete legend.

| Data | nD | eO | eF | mV | eT | eS | eP | e$\Sigma$ | eR | eH | eA | eB | nC | nE | nG |
|---|---|---|---|---|---|---|---|---|---|---|---|---|---|---|---|
| cub06 | **1.00** | **1.00** | **1.00** | **1.00** | **1.00** | **1.00** | **1.00** | **1.00** | **1.00** | **1.00** | **1.00** | **1.00** | **1.00** | **1.00** | **1.00** |
| cub08 | 0.90 | 0.90 | 0.90 | 0.90 | 0.90 | **1.00** | **1.00** | **1.00** | 0.90 | 0.90 | 0.90 | **1.00** | 0.90 | 0.90 | **1.00** |
| cub10 | 0.81 | 0.99 | 0.88 | 0.85 | 0.88 | **1.00** | 0.98 | **1.00** | 0.81 | 0.81 | 0.85 | **1.00** | 0.88 | 0.88 | 0.99 |
| cub12 | 0.80 | **1.00** | 0.87 | 0.87 | 0.87 | 0.99 | 0.99 | **1.00** | 0.80 | 0.80 | 0.84 | **1.00** | 0.87 | 0.87 | **1.00** |
| cub14 | 0.79 | **1.00** | 0.86 | 0.89 | 0.86 | 0.99 | **1.00** | **1.00** | 0.79 | 0.79 | 0.83 | **1.00** | 0.86 | 0.86 | **1.00** |

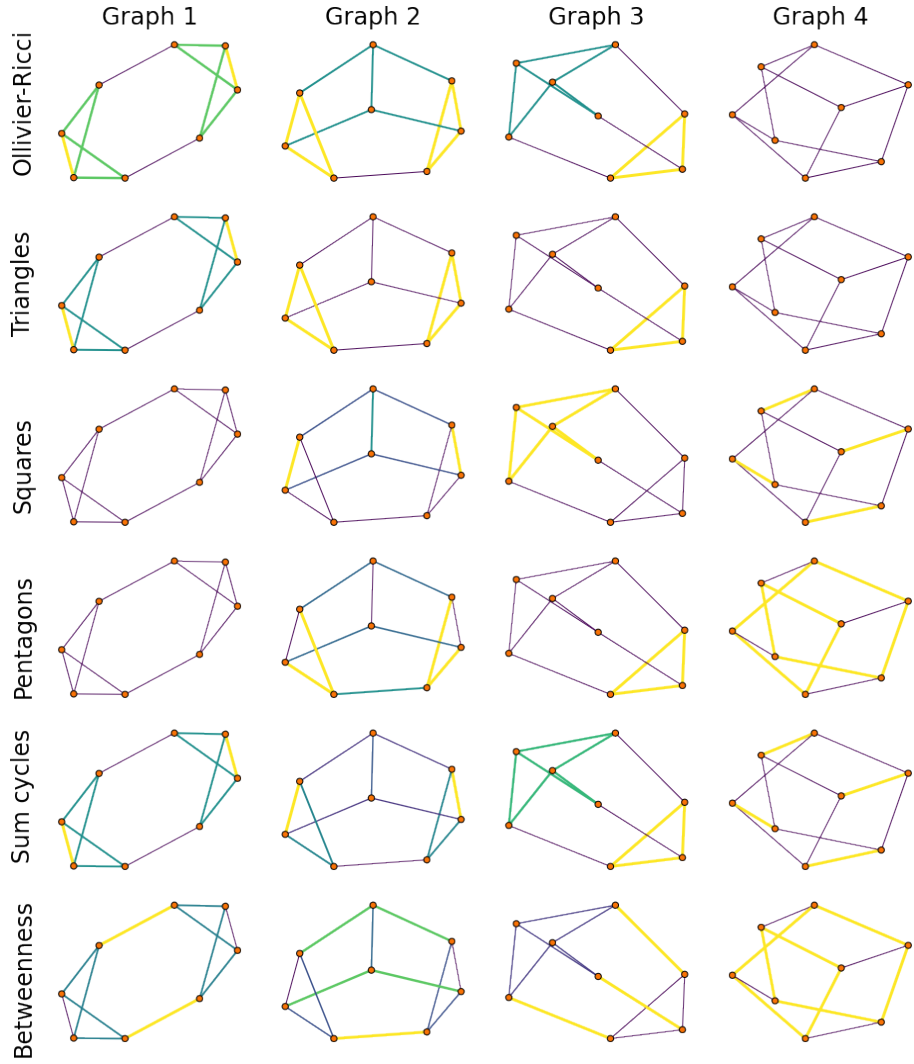

FIG. S3: Four different graphs from the *cubic* dataset. Each row corresponds to a different measure of edge weights. The lighter the color, the larger the values are. For each edge weighting, the weights are normalized to improve the visualization. We see that this dataset is distinguishable by all weighting filtrations.

TABLE S5: Success rate for distinguishing pairs of non-isomorphic *minimal Cayley* graphs when using different filtrations at expansion level of the graph $k = 1$. Bold indicates the best performing algorithm. Filtrations are denoted by their abbreviations; see Table I for the complete legend.

| Data | nD | eO | eF | mV | eT | eS | eP | eΣ | eR | eH | eA | eB | nC | nE | nG |
|---|---|---|---|---|---|---|---|---|---|---|---|---|---|---|---|
| cay12 | 0.67 | 0.71 | 0.95 | 0.00 | 0.71 | 0.95 | 0.29 | **1.00** | 0.67 | 0.67 | 0.95 | **1.00** | 0.67 | 0.86 | 0.95 |
| cay16 | 0.83 | 0.42 | 0.83 | 0.00 | 0.00 | 0.97 | 0.00 | 0.97 | 0.83 | 0.83 | 0.83 | **1.00** | 0.00 | 0.83 | **1.00** |
| cay20 | 0.61 | 0.46 | 0.61 | 0.00 | 0.00 | **1.00** | 0.46 | **1.00** | 0.61 | 0.61 | 0.61 | **1.00** | 0.00 | 0.61 | 0.96 |
| cay24 | 0.65 | 0.82 | 0.86 | 0.00 | 0.54 | 0.93 | 0.35 | 0.98 | 0.65 | 0.65 | 0.83 | **1.00** | 0.52 | 0.74 | 0.97 |
| cay32 | 0.76 | 0.81 | 0.76 | 0.00 | 0.00 | 0.98 | 0.00 | 0.98 | 0.76 | 0.76 | 0.76 | **1.00** | 0.00 | 0.76 | 0.98 |
| cay36 | 0.69 | 0.87 | 0.84 | 0.00 | 0.58 | 0.95 | 0.34 | 0.98 | 0.69 | 0.69 | 0.82 | **1.00** | 0.57 | 0.75 | 0.97 |
| cay60 | 0.69 | 0.90 | 0.78 | 0.00 | 0.28 | 0.91 | 0.50 | 0.96 | 0.69 | 0.69 | 0.77 | **1.00** | 0.28 | 0.73 | 0.95 |
| cay63 | 0.49 | 0.89 | 0.73 | 0.00 | 0.67 | 0.78 | 0.35 | 0.90 | 0.49 | 0.49 | 0.67 | **1.00** | 0.67 | 0.73 | 0.90 |

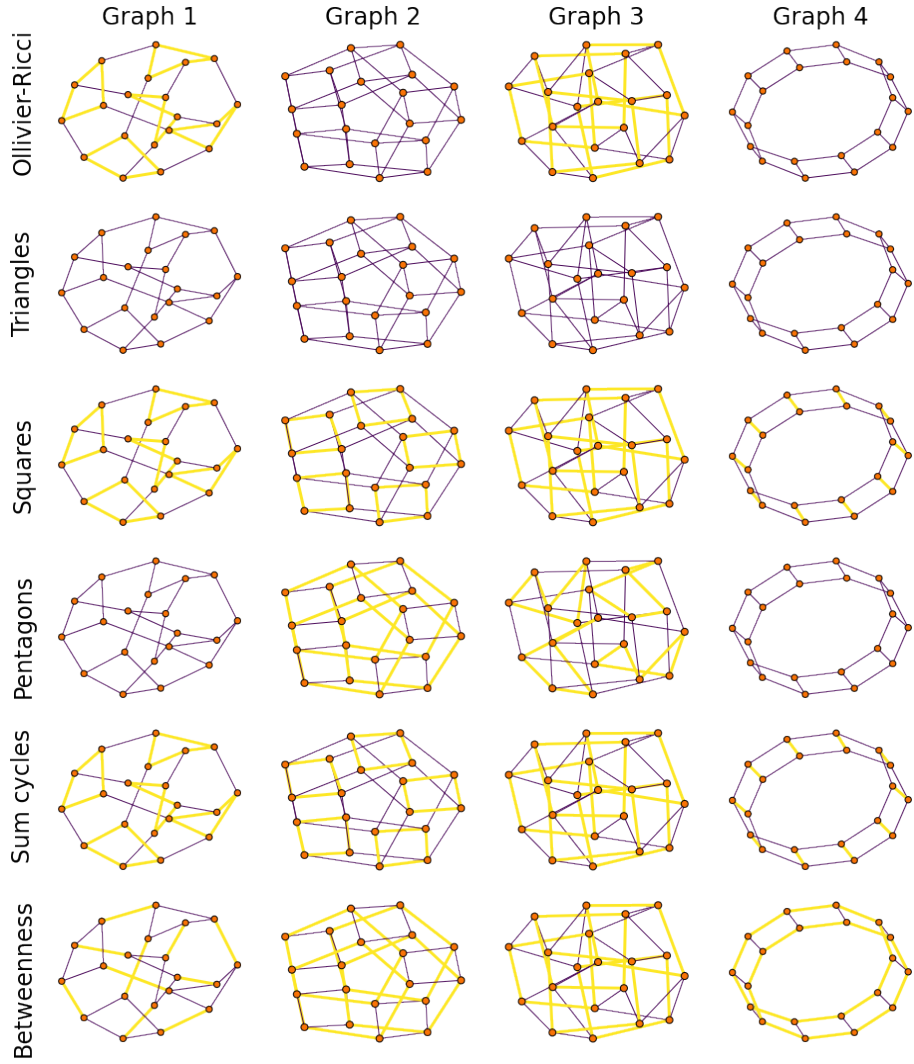

FIG. S4: Four different graphs from the *minimal Cayley* dataset. Each row corresponds to a different measure of edge weights. The lighter the color, the larger the values are. For each edge weighting, the weights are normalized to improve the visualization. Ollivier–Ricci does not distinguish between graphs 2 and 4. Similarly, triangle weighting provides no differences between graphs, and pentagons do not distinguish between graphs 1 and 4. Instead, betweenness and the sum of cycle densities allow us to distinguish all of them.

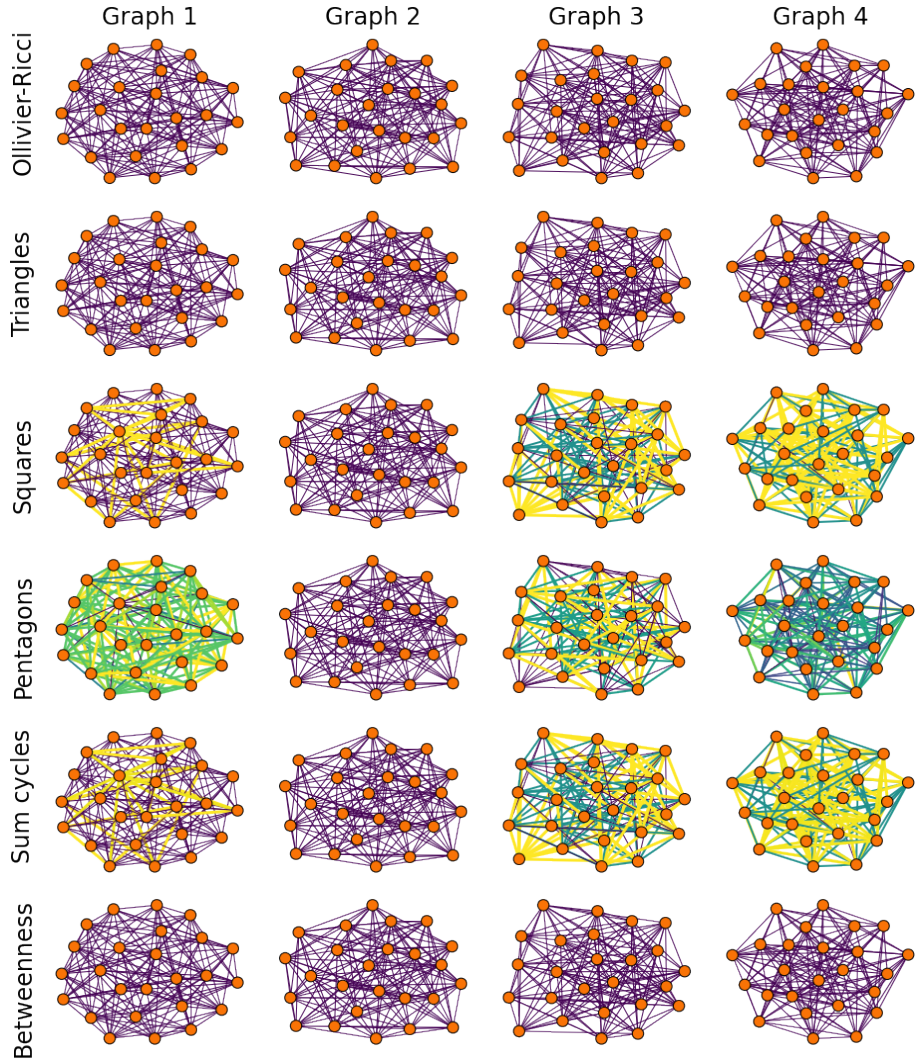

FIG. S5: Four different graphs from the *strongly regular* dataset. Each row corresponds to a different measure of edge weights. The lighter the color, the larger the values are. For each edge weighting, the weights are normalized to improve the visualization. While Ollivier–Ricci curvature, triangles, and betweenness weightings provide no differences between graphs, cycle densities allow us to distinguish all of them.

TABLE S6: Success rate for distinguishing pairs of non-isomorphic *minimal Cayley* graphs when using different filtrations at expansion level of the graph $k = 2$. Bold indicates the best performing algorithm. Filtrations are denoted by their abbreviations; see Table I for the complete legend.

| Data | nD | eO | eF | mV | eT | eS | eP | eΣ | eR | eH | eA | eB | nC | nE | nG |
|---|---|---|---|---|---|---|---|---|---|---|---|---|---|---|---|
| cay12 | 0.95 | 0.95 | 0.95 | 0.95 | 0.95 | **1.00** | **1.00** | **1.00** | 0.95 | 0.95 | 0.95 | **1.00** | 0.95 | 0.95 | **1.00** |
| cay16 | 0.83 | 0.92 | 0.83 | 0.94 | 0.83 | **1.00** | 0.83 | **1.00** | 0.83 | 0.83 | 0.83 | **1.00** | 0.83 | 0.83 | **1.00** |
| cay20 | 0.61 | 0.79 | 0.61 | 0.89 | 0.61 | **1.00** | 0.86 | **1.00** | 0.61 | 0.61 | 0.61 | **1.00** | 0.61 | 0.61 | 0.96 |
| cay24 | 0.83 | 0.93 | 0.86 | 0.93 | 0.86 | 0.97 | 0.90 | 0.98 | 0.83 | 0.83 | 0.83 | **1.00** | 0.83 | 0.83 | 0.97 |
| cay32 | 0.76 | 0.94 | 0.76 | 0.90 | 0.76 | 0.99 | 0.76 | 0.99 | 0.76 | 0.76 | 0.76 | **1.00** | 0.76 | 0.76 | 0.98 |
| cay36 | 0.84 | 0.95 | 0.84 | 0.94 | 0.84 | 0.98 | 0.89 | 0.99 | 0.84 | 0.84 | 0.82 | **1.00** | 0.84 | 0.84 | 0.97 |
| cay60 | 0.77 | 0.95 | 0.78 | 0.97 | 0.78 | 0.96 | 0.87 | 0.97 | 0.77 | 0.77 | 0.77 | **1.00** | 0.77 | 0.77 | 0.95 |
| cay63 | 0.73 | 0.93 | 0.73 | 0.96 | 0.73 | 0.90 | 0.76 | 0.93 | 0.73 | 0.73 | 0.67 | **1.00** | 0.73 | 0.73 | 0.90 |

TABLE S7: Success rate for distinguishing pairs of non-isomorphic *strongly regular* graphs when using different filtrations at expansion level of the graph $k = 1$. Bold indicates the best-performing algorithm. Filtrations are denoted by their abbreviations; see Table I for the complete legend.

| Data | nD | eO | eF | mV | eT | eS | eP | e$\Sigma$ | eR | eH | eA | eB | nC | nE | nG |
|---|---|---|---|---|---|---|---|---|---|---|---|---|---|---|---|
| sr16622  | 0.00 | **1.00** | 0.00 | 0.00 | 0.00 | **1.00** | **1.00** | **1.00** | 0.00 | 0.00 | 0.00 | 0.00 | 0.00 | 0.00 | **1.00** |
| sr251256 | 0.00 | 0.00 | 0.00 | 0.00 | 0.00 | 0.94 | **1.00** | **1.00** | 0.00 | 0.00 | 0.00 | 0.00 | 0.00 | 0.00 | 0.65 |
| sr261034 | 0.00 | 0.00 | 0.00 | 0.00 | 0.00 | 0.00 | **1.00** | **1.00** | 0.00 | 0.00 | 0.00 | 0.00 | 0.00 | 0.00 | 0.84 |
| sr281264 | 0.00 | 0.83 | 0.00 | 0.00 | 0.00 | 0.83 | **1.00** | 0.83 | 0.00 | 0.00 | 0.00 | 0.00 | 0.00 | 0.00 | **1.00** |
| sr291467 | 0.00 | 0.00 | 0.00 | 0.00 | 0.00 | 0.94 | **1.00** | **1.00** | 0.00 | 0.00 | 0.00 | 0.00 | 0.00 | 0.00 | 0.80 |
| sr351668 | 0.00 | 0.00 | 0.00 | 0.00 | 0.00 | **1.00** | **1.00** | **1.00** | 0.00 | 0.00 | 0.00 | 0.00 | 0.00 | 0.00 | 0.98 |
| sr351899 | 0.00 | 0.00 | 0.00 | 0.00 | 0.00 | 0.98 | **1.00** | **1.00** | 0.00 | 0.00 | 0.00 | 0.00 | 0.00 | 0.00 | 0.76 |
| sr361446 | 0.00 | 0.00 | 0.00 | 0.00 | 0.00 | 0.12 | **1.00** | **1.00** | 0.00 | 0.00 | 0.00 | 0.00 | 0.00 | 0.00 | 0.84 |
| sr401224 | 0.00 | 0.00 | 0.00 | 0.00 | 0.00 | 0.89 | 0.97 | **0.99** | 0.00 | 0.00 | 0.00 | 0.00 | 0.00 | 0.00 | 0.43 |

TABLE S8: Success rate for distinguishing pairs of non-isomorphic *strongly regular* graphs when using different filtrations at expansion level of the graph $k = 2$. Bold indicates the best-performing algorithm. Filtrations are denoted by their abbreviations; see Table I for the complete legend.

| Data | nD | eO | eF | mV | eT | eS | eP | e$\Sigma$ | eR | eH | eA | eB | nC | nE | nG |
|---|---|---|---|---|---|---|---|---|---|---|---|---|---|---|---|
| sr16622  | **1.00** | **1.00** | **1.00** | **1.00** | **1.00** | **1.00** | **1.00** | **1.00** | **1.00** | **1.00** | **1.00** | **1.00** | **1.00** | **1.00** | **1.00** |
| sr251256 | 0.00 | 0.00 | 0.00 | 0.00 | 0.00 | **1.00** | **1.00** | **1.00** | 0.00 | 0.00 | 0.00 | 0.90 | 0.00 | 0.00 | 0.95 |
| sr261034 | 0.20 | 0.20 | 0.20 | 0.20 | 0.20 | 0.76 | **1.00** | **1.00** | 0.20 | 0.20 | 0.20 | 0.93 | 0.20 | 0.20 | 0.96 |
| sr281264 | 0.00 | **1.00** | 0.00 | 0.00 | 0.00 | **1.00** | **1.00** | **1.00** | 0.00 | 0.00 | 0.00 | **1.00** | 0.00 | 0.00 | **1.00** |
| sr291467 | 0.00 | 0.00 | 0.00 | 0.00 | 0.00 | **1.00** | **1.00** | **1.00** | 0.00 | 0.00 | 0.00 | 0.77 | 0.00 | 0.00 | 0.94 |
| sr351668 | 0.00 | 0.00 | 0.00 | 0.00 | 0.00 | **1.00** | **1.00** | **1.00** | 0.00 | 0.00 | 0.00 | 0.00 | 0.00 | 0.00 | **1.00** |
| sr351899 | 0.00 | 0.00 | 0.00 | 0.00 | 0.00 | **1.00** | **1.00** | **1.00** | 0.00 | 0.00 | 0.00 | 0.81 | 0.00 | 0.00 | 0.81 |
| sr361446 | 0.02 | 0.02 | 0.02 | 0.02 | 0.02 | 0.88 | **1.00** | **1.00** | 0.02 | 0.02 | 0.02 | 0.92 | 0.02 | 0.02 | 0.98 |
| sr401224 | 0.93 | 0.93 | 0.93 | 0.93 | 0.93 | 0.98 | **1.00** | **1.00** | 0.93 | 0.93 | 0.93 | 0.94 | 0.93 | 0.93 | 0.99 |

TABLE S9: Success rate for distinguishing pairs of instances of the BREC data set when using different filtrations at expansion level of the graph $k = 1$. Legend and number of graphs per category: B (Basic, 60), R (Regular, 100), E (Extension, 100), C (CFI, 100), 4, (4-Vertex Condition, 20), D (Distance-Regular, 20) graphs, respectively and A (average over full data set, 400 graphs). Bold indicates the best-performing algorithm. Filtrations are denoted by their abbreviations; see Table I for the complete legend.

| Data | Filtration ($k = 1$) | | | | | | | | | | | | | | |
|---|---|---|---|---|---|---|---|---|---|---|---|---|---|---|---|
|  | nD | eO | eF | mV | eT | eS | eP | e$\Sigma$ | eR | eH | eA | eB | nC | nE | nG |
| B | 0.033 | 0.933 | 0.750 | 0.000 | 0.933 | 0.867 | 0.533 | 0.983 | 0.000 | 0.000 | 0.983 | **1.000** | 0.667 | 0.800 | **1.000** |
| R | 0.000 | 0.420 | 0.320 | 0.000 | 0.360 | 0.900 | 0.920 | **1.000** | 0.000 | 0.000 | 0.410 | 0.500 | 0.410 | 0.410 | 0.860 |
| E | 0.070 | 0.760 | 0.440 | 0.000 | 0.510 | 0.690 | 0.490 | 0.950 | 0.000 | 0.000 | 0.670 | **1.000** | 0.000 | 0.100 | 0.850 |
| C | 0.030 | 0.030 | 0.030 | 0.030 | 0.030 | 0.030 | 0.030 | 0.030 | 0.030 | 0.030 | 0.030 | **0.100** | 0.030 | 0.030 | 0.030 |
| 4 | 0.000 | 0.000 | 0.000 | 0.000 | 0.000 | 0.000 | **0.600** | **0.600** | 0.000 | 0.000 | 0.000 | 0.000 | 0.000 | 0.000 | 0.000 |
| D | 0.000 | 0.000 | 0.000 | 0.000 | 0.000 | 0.000 | 0.000 | 0.000 | 0.000 | 0.000 | 0.000 | 0.000 | 0.000 | 0.000 | 0.000 |
| A | 0.030 | 0.442 | 0.310 | 0.008 | 0.365 | 0.535 | 0.470 | **0.672** | 0.008 | 0.008 | 0.425 | 0.550 | 0.210 | 0.255 | 0.585 |

TABLE S10: Success rate for distinguishing pairs of instances of the BREC data set when using different filtrations at expansion level of the graph $k = 2$. Legend and number of graphs per category: B (Basic, 60), R (Regular, 100), E (Extension, 100), C (CFI, 100), 4 (4-Vertex Condition, 20), D (Distance-Regular, 20) graphs, respectively and A (average over full data set, 400 graphs). Bold indicates the best-performing algorithm. Filtrations are denoted by their abbreviations; see Table I for the complete legend.

| Data | Filtration ($k = 2$) | | | | | | | | | | | | | |
|---|---|---|---|---|---|---|---|---|---|---|---|---|---|---|
| | nD | eO | eF | mV | eT | eS | eP | eΣ | eR | eH | eA | eB | nC | nE | nG |
| B | 0.783 | **1.000** | 0.983 | 0.517 | 0.983 | 0.867 | 0.733 | **1.000** | 0.517 | 0.517 | 0.983 | **1.000** | 0.750 | 0.883 | **1.000** |
| R | 0.390 | 0.540 | 0.500 | 0.390 | 0.500 | 0.980 | 0.970 | **1.000** | 0.390 | 0.390 | 0.500 | 0.520 | 0.510 | 0.510 | 0.960 |
| E | 0.260 | 0.920 | 0.590 | 0.110 | 0.640 | 0.760 | 0.580 | 0.980 | 0.120 | 0.120 | 0.700 | **1.000** | 0.120 | 0.250 | 0.890 |
| C | 0.030 | 0.030 | 0.030 | 0.030 | 0.030 | 0.030 | 0.030 | 0.030 | 0.030 | 0.030 | 0.030 | **0.100** | 0.030 | 0.030 | 0.030 |
| 4 | 0.000 | 0.000 | 0.000 | 0.000 | 0.000 | 0.000 | **0.800** | **0.800** | 0.000 | 0.000 | 0.000 | 0.000 | 0.000 | 0.000 | 0.000 |
| D | 0.000 | 0.000 | 0.000 | 0.000 | 0.000 | 0.000 | 0.000 | 0.000 | 0.000 | 0.000 | 0.000 | 0.000 | 0.000 | 0.000 | 0.000 |
| A | 0.290 | 0.523 | 0.427 | 0.210 | 0.440 | 0.573 | 0.545 | **0.693** | 0.213 | 0.213 | 0.455 | 0.555 | 0.278 | 0.330 | 0.620 |

TABLE S11: Success rate for distinguishing pairs of instances of the BREC data set when using different filtrations at expansion level of the graph $k = 3$. Legend and number of graphs per category: B (Basic, 60), R (Regular, 100), E (Extension, 100), C (CFI, 100), 4 (4-Vertex Condition, 20), D (Distance-Regular, 20) graphs, respectively and A (average over full data set, 400 graphs). Bold indicates the best-performing algorithm. Filtrations are denoted by their abbreviations; see Table I for the complete legend.

| Data | Filtration ($k = 3$) | | | | | | | | | | | | | |
|---|---|---|---|---|---|---|---|---|---|---|---|---|---|---|
| | nD | eO | eF | mV | eT | eS | eP | eΣ | eR | eH | eA | eB | nC | nE | nG |
| B | 0.833 | **1.000** | 0.983 | 0.583 | 0.983 | 0.867 | 0.800 | **1.000** | 0.583 | 0.583 | 0.983 | **1.000** | 0.78 | 0.917 | **1.000** |
| R | 0.850 | 0.930 | 0.910 | 0.850 | 0.910 | **1.000** | **1.000** | **1.000** | 0.850 | 0.850 | 0.910 | 0.930 | 0.920 | 0.920 | **1.000** |
| E | 0.290 | 0.920 | 0.590 | 0.160 | 0.640 | 0.760 | 0.710 | 0.980 | 0.150 | 0.150 | 0.700 | **1.000** | 0.150 | 0.280 | 0.890 |
| C | 0.030 | 0.030 | 0.030 | 0.030 | 0.030 | 0.030 | 0.030 | 0.030 | 0.030 | 0.030 | 0.030 | **0.100** | 0.030 | 0.030 | 0.030 |
| 4 | **1.000** | **1.000** | **1.000** | **1.000** | **1.000** | **1.000** | **1.000** | **1.000** | **1.000** | **1.000** | **1.000** | **1.000** | **1.000** | **1.000** | **1.000** |
| D | 0.000 | 0.000 | 0.000 | **0.050** | 0.000 | 0.000 | 0.000 | 0.000 | 0.000 | 0.000 | 0.000 | 0.000 | 0.000 | 0.000 | 0.000 |
| A | 0.467 | 0.670 | 0.580 | 0.400 | 0.592 | 0.628 | 0.605 | 0.703 | 0.395 | 0.395 | 0.607 | **0.708** | 0.442 | 0.495 | 0.680 |

TABLE S12: Success rate for distinguishing pairs of instances of the BREC data set when using different filtrations at expansion level of the graph $k = 4$. Like in [3], due to combinatorial constraints, we did not calculate the Vietoris–Rips filtration for $k = 4$. Legend and number of graphs per category: B (Basic, 60), R (Regular, 100), E (Extension, 100), C (CFI, 100), 4 (4-Vertex Condition, 20), D (Distance-Regular, 20) graphs, respectively and A (average over full data set, 400 graphs). Bold indicates the best-performing algorithm. Filtrations are denoted by their abbreviations; see Table I for the complete legend.

| Data | Filtration ($k = 4$) | | | | | | | | | | | | | |
|---|---|---|---|---|---|---|---|---|---|---|---|---|---|---|
| | nD | eO | eF | mV | eT | eS | eP | eΣ | eR | eH | eA | eB | nC | nE | nG |
| B | 0.833 | **1.000** | 0.983 | - | 0.983 | 0.867 | 0.733 | **1.000** | 0.583 | 0.583 | 0.983 | **1.000** | 0.783 | 0.917 | **1.000** |
| R | 0.890 | 0.970 | 0.950 | - | 0.950 | **1.000** | **1.000** | **1.000** | 0.890 | 0.890 | 0.950 | 0.970 | 0.960 | 0.960 | **1.000** |
| E | 0.290 | 0.920 | 0.590 | - | 0.640 | 0.760 | 0.580 | 0.980 | 0.150 | 0.150 | 0.700 | **1.000** | 0.150 | 0.280 | 0.890 |
| C | 0.030 | 0.030 | 0.030 | - | 0.030 | 0.030 | 0.030 | 0.030 | 0.030 | 0.030 | 0.030 | **0.100** | 0.030 | 0.030 | 0.030 |
| 4 | **1.000** | **1.000** | **1.000** | - | **1.000** | **1.000** | **1.000** | **1.000** | **1.000** | **1.000** | **1.000** | **1.000** | **1.000** | **1.000** | **1.000** |
| D | 0.000 | 0.000 | 0.000 | - | 0.000 | 0.000 | 0.000 | 0.000 | 0.000 | 0.000 | 0.000 | 0.000 | 0.000 | 0.000 | 0.000 |
| A | 0.479 | 0.680 | 0.588 | - | 0.603 | 0.627 | 0.563 | 0.703 | 0.405 | 0.405 | 0.618 | **0.718** | 0.453 | 0.505 | 0.680 |

TABLE S13: Summary statistics for all graph properties in the ogbg-molhiv dataset [2].

|                             | Minimum | Maximum  | Mean    | Variance |
|-----------------------------|---------|----------|---------|----------|
| radius                      | 1.0000  | 47.0000  | 6.3046  | 6.9007   |
| girth                       | 3.0000  | 105.0000 | 5.9607  | 10.9867  |
| diameter                    | 1.0000  | 93.0000  | 11.9847 | 26.7765  |
| avg. clustering coefficient | 0.0000  | 0.6364   | 0.0020  | 0.0003   |
| avg. shortest path length   | 1.0000  | 35.1918  | 4.9761  | 3.3375   |
| avg. degree                 | 1.0000  | 3.5200   | 2.1446  | 0.0131   |
| minimum degree              | 0.0000  | 2.0000   | 0.9660  | 0.0754   |
| maximum degree              | 1.0000  | 10.0000  | 3.3817  | 0.3015   |
| avg. betweenness centrality | 0.0000  | 0.3333   | 0.1740  | 0.0012   |
| avg. closeness centrality   | 0.0295  | 1.0000   | 0.2270  | 0.0048   |
| degree heterogeneity        | 0.0000  | 1.3333   | 0.1291  | 0.0037   |

TABLE S14: Property prediction accuracy for dimension $k = 1$, when predicting the following properties of graphs: radius, diameter and girth, using the ogbg-molhiv data set [2] and the filtrations studied in this paper. We exclude mV results due to out-of-resources errors during execution. Filtrations are denoted by their abbreviations; see Table I for the complete legend.

| Data     | $k = 1$ |       |       |    |       |       |       |       |       |       |       |       |       |       |       |
|----------|-------|-------|-------|----|-------|-------|-------|-------|-------|-------|-------|-------|-------|-------|-------|
|          | nD    | eO    | eF    | mV | eT    | eS    | eP    | e$\Sigma$ | eR    | eH    | eA    | eB    | nC    | nE    | nG    |
| radius   | 0.163 | 0.211 | 0.057 | -  | 0.160 | 0.160 | 0.156 | 0.159 | 0.046 | 0.034 | 0.008 | 0.014 | 0.130 | 0.099 | 0.145 |
| diameter | 0.089 | 0.110 | 0.049 | -  | 0.077 | 0.077 | 0.077 | 0.078 | 0.030 | 0.018 | 0.007 | 0.009 | 0.055 | 0.044 | 0.072 |
| girth    | 0.001 | 0.113 | 0.340 | -  | 0.481 | 0.483 | 0.447 | 0.436 | 0.334 | 0.506 | 0.037 | 0.247 | 0.476 | 0.459 | 0.458 |

TABLE S15: Mean absolute error for dimension $k = 1$ when predicting the following properties of graphs: average closeness centrality, average degree, average clustering coefficient, degree heterogeneity, average betweenness centrality, minimum degree, maximum degree, and average shortest path length, using the ogbg-molhiv data set [2] and the filtrations studied in this paper. We exclude mV results due to out-of-resources errors during execution. Filtrations are denoted by their abbreviations; see Table I for the complete legend.

| Data | $k = 1$ | | | | | | | | | | | | | | |
|------|----|----|----|----|----|----|----|----|----|----|----|----|----|----|----|
|      | nD | eO | eF | mV | eT | eS | eP | e$\Sigma$ | eR | eH | eA | eB | nC | nE | nG |
| avg. closeness centrality   | 0.154 | 0.129 | 0.084 | - | 0.085 | 0.085 | 0.085 | 0.085 | 0.081 | 0.317 | 0.524 | 0.098 | 0.089 | 0.079 | 0.083 |
| avg. degree                 | 0.448 | 0.181 | 0.102 | - | 0.098 | 0.098 | 0.095 | 0.094 | 0.183 | 0.186 | 0.977 | 0.142 | 0.145 | 0.089 | 0.110 |
| avg. clustering coefficient | 0.006 | 0.093 | 0.011 | - | 0.956 | 0.956 | 0.956 | 0.956 | 0.079 | 0.111 | 0.012 | 0.167 | 0.008 | 0.125 | 0.010 |
| degree heterogeneity        | 0.065 | 0.114 | 0.070 | - | 0.045 | 0.044 | 0.046 | 0.047 | 0.047 | 0.040 | 0.648 | 0.045 | 0.083 | 0.050 | 0.056 |
| avg. betweenness centrality | 0.136 | 0.040 | 0.034 | - | 0.029 | 0.029 | 0.030 | 0.030 | 0.065 | 0.052 | 0.076 | 0.144 | 0.041 | 0.038 | 0.036 |
| minumum degree              | 0.133 | 0.579 | 0.120 | - | 0.135 | 0.133 | 0.135 | 0.135 | 0.264 | 0.523 | 0.939 | 0.708 | 0.621 | 0.339 | 0.138 |
| maximum degree              | 1.322 | 1.431 | 0.610 | - | 0.501 | 0.504 | 0.501 | 0.500 | 0.469 | 0.485 | 3.206 | 0.574 | 0.586 | 0.600 | 0.549 |
| avg. shortest path length   | 1.606 | 1.314 | 2.698 | - | 1.798 | 1.798 | 1.785 | 1.800 | 2.099 | 2.613 | 3.337 | 8.501 | 1.755 | 1.917 | 1.959 |

TABLE S16: Property prediction accuracy for dimension $k = 2$, when predicting the following properties of graphs: radius, diameter and girth, using the ogbg-molhiv dataset [2] and the filtrations studied in this paper. We exclude mV results due to out-of-resources errors during execution. Filtrations are denoted by their abbreviations; see Table I for the complete legend.

| Data | $k = 2$ | | | | | | | | | | | | | | |
|---|---|---|---|---|---|---|---|---|---|---|---|---|---|---|---|
| | nD | eO | eF | mV | eT | eS | eP | e$\Sigma$ | eR | eH | eA | eB | nC | nE | nG |
| radius | 0.200 | 0.182 | 0.097 | - | 0.160 | 0.159 | 0.155 | 0.156 | 0.050 | 0.036 | 0.008 | 0.014 | 0.083 | 0.090 | 0.138 |
| diameter | 0.103 | 0.075 | 0.060 | - | 0.077 | 0.076 | 0.079 | 0.078 | 0.033 | 0.016 | 0.007 | 0.006 | 0.046 | 0.045 | 0.073 |
| girth | 0.144 | 0.212 | 0.333 | - | 0.480 | 0.488 | 0.451 | 0.437 | 0.238 | 0.475 | 0.018 | 0.510 | 0.031 | 0.381 | 0.475 |

TABLE S17: Mean absolute error for dimension $k = 2$ when predicting the following properties of graphs: average closeness centrality, average degree, average clustering coefficient, degree heterogeneity, average betweenness centrality, minimum degree, maximum degree, and average shortest path length, using the ogbg-molhiv data set [2] and the filtrations studied in this paper. We exclude mV results due to out-of-resources errors during execution.

| Data | $k = 2$ | | | | | | | | | | | | | | |
|---|---|---|---|---|---|---|---|---|---|---|---|---|---|---|---|
| | nD | eO | eF | mV | eT | eS | eP | e$\Sigma$ | eR | eH | eA | eB | nC | nE | nG |
| avg. closeness centrality | 0.126 | 0.129 | 0.087 | - | 0.085 | 0.085 | 0.085 | 0.085 | 0.079 | 0.310 | 0.527 | 0.102 | 0.096 | 0.079 | 0.082 |
| avg. degree | 0.142 | 0.136 | 0.109 | - | 0.098 | 0.098 | 0.096 | 0.095 | 0.224 | 0.321 | 0.987 | 0.136 | 0.077 | 0.088 | 0.100 |
| avg. clustering coefficient | 0.030 | 0.150 | 0.057 | - | 0.956 | 0.956 | 0.956 | 0.956 | 0.018 | 0.011 | 0.006 | 0.131 | 0.056 | 0.132 | 0.956 |
| degree heterogeneity | 0.050 | 0.119 | 0.103 | - | 0.045 | 0.044 | 0.046 | 0.047 | 0.057 | 0.040 | 0.649 | 0.049 | 0.063 | 0.050 | 0.056 |
| avg. betweenness centrality | 0.034 | 0.037 | 0.049 | - | 0.029 | 0.029 | 0.030 | 0.030 | 0.066 | 0.054 | 0.079 | 0.143 | 0.051 | 0.036 | 0.035 |
| minimum degree | 0.133 | 0.525 | 0.095 | - | 0.134 | 0.131 | 0.133 | 0.139 | 0.281 | 0.524 | 0.943 | 0.741 | 0.153 | 0.300 | 0.142 |
| maximum degree | 1.217 | 1.844 | 1.839 | - | 0.504 | 0.498 | 0.500 | 0.503 | 0.483 | 0.492 | 2.448 | 1.628 | 0.649 | 0.582 | 0.548 |
| avg. shortest path length | 1.617 | 1.312 | 2.478 | - | 1.799 | 1.798 | 1.788 | 1.800 | 2.108 | 2.580 | 3.329 | 8.712 | 2.018 | 1.909 | 1.952 |

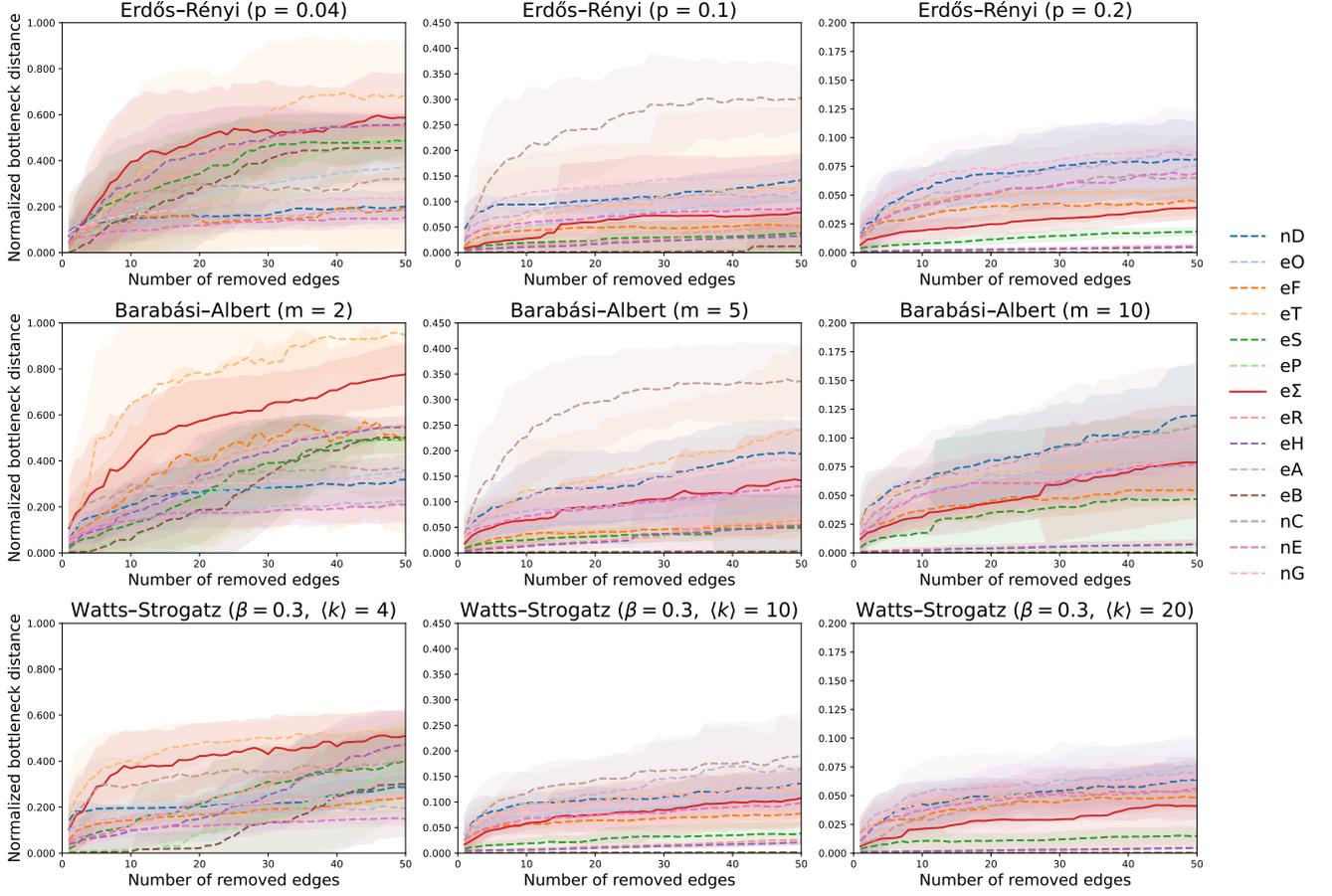

FIG. S6: Sensitivity of network filtrations across random graph models under edge removal perturbations. Each panel shows the average normalized bottleneck distance in homology dimension $k = 0$, computed over 50 runs of independently generated graphs with $n = 100$, as a function of the number of removed edges. Results are presented for ER graphs with connection probabilities $p \in \{0.04, 0.1, 0.2\}$, BA graphs with attachment parameters $m \in \{2, 5, 10\}$, and WS small-world graphs with rewiring probability $\beta = 0.3$ and neighborhood sizes $\langle s \rangle \in \{4, 10, 20\}$. All configurations yield comparable average degrees across models. Several graph-based filtrations were evaluated, including degree, curvature-based, and motif-density measures. The results reveal that network sensitivity to edge removal decreases with increasing graph density and varies substantially across filtration types and generative models. Filtrations are denoted by their abbreviations; see Table I for the complete legend.

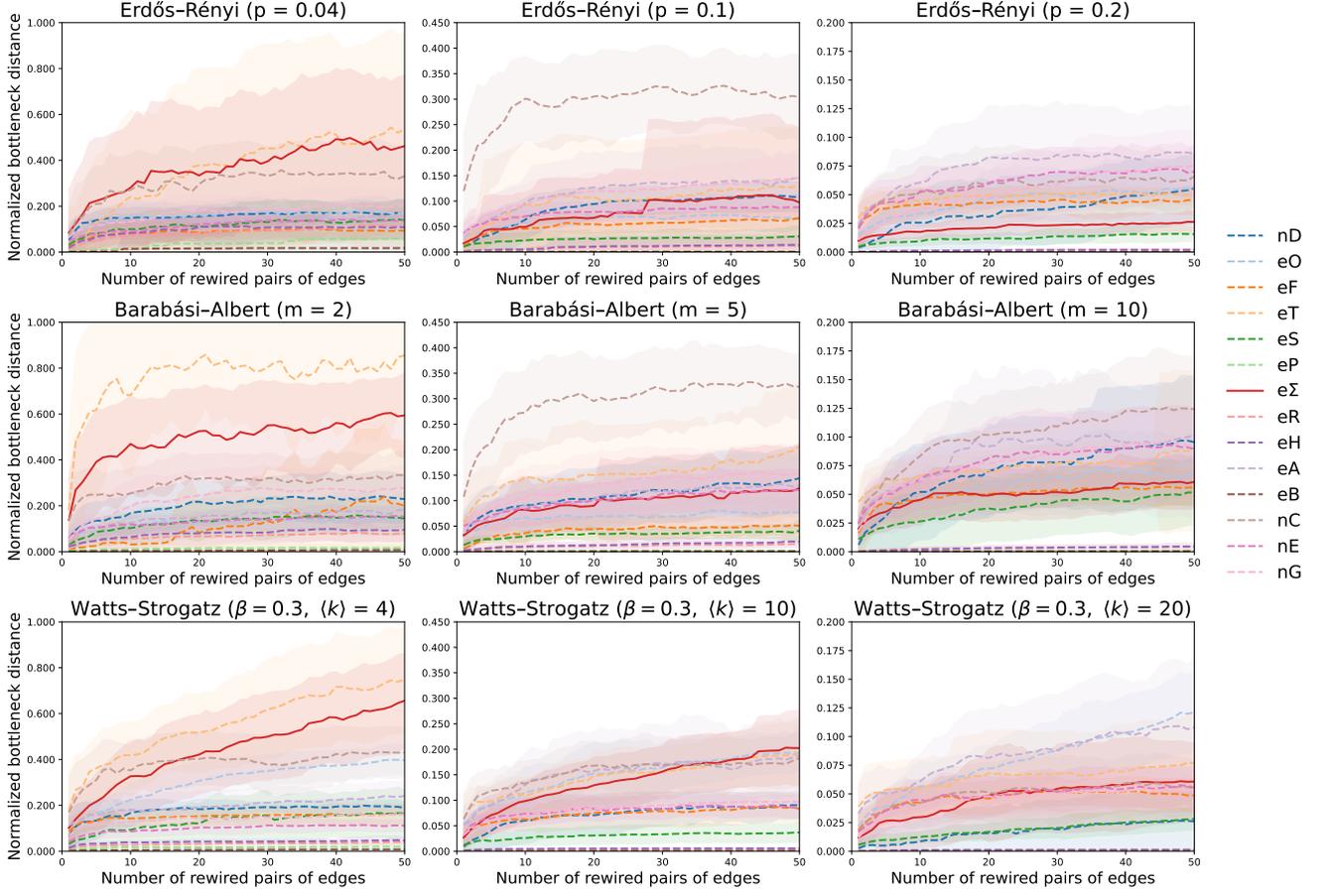

FIG. S7: Sensitivity of network filtrations across random graph models under edge rewiring perturbations. Each panel shows the average normalized bottleneck distance in homology dimension $k = 0$, computed over 50 runs of independently generated graphs with $n = 100$, as a function of the number of removed edges. Results are presented for ER graphs with connection probabilities $p \in \{0.04, 0.1, 0.2\}$, BA graphs with attachment parameters $m \in \{2, 5, 10\}$, and WS small-world graphs with rewiring probability $\beta = 0.3$ and neighborhood sizes $\langle s \rangle \in \{4, 10, 20\}$. All configurations yield comparable average degrees across models. Several graph-based filtrations were evaluated, including degree, curvature-based, and motif-density measures. The results reveal that network sensitivity to edge rewiring decreases with increasing graph density and varies substantially across filtration types and generative models. Filtrations are denoted by their abbreviations; see Table I for the complete legend.

## SUPPLEMENTARY REFERENCES



\end{document}